\documentclass[twoside,a4paper,11pt]{amsart}
\usepackage{amsmath,amssymb,amsthm,amscd}
\usepackage[]{latexsym,amssymb,amsmath,amsfonts, amsthm, verbatim}
\usepackage[all,cmtip,ps]{xy}
\usepackage[latin1]{inputenc}
\usepackage[english]{babel}

\usepackage{graphicx}

\newcommand{\add}{\mathrm{add}}
\newcommand{\Add}{\mathrm{Add}}

\newcommand{\Gen}{\mathrm{Gen}}

\newcommand{\rmod}{\mathrm{Mod-}}

\newcommand{\M}{\mbox{\rm Mod-$R$}}
\newcommand{\m}{\mbox{\rm mod-$R$}}

\newcommand{\Z}{\mathbb{Z}}

\newcommand{\N}{\mathbb{N}}

\DeclareMathOperator{\End}{End}
\DeclareMathOperator{\Ker}{Ker}
\DeclareMathOperator{\Img}{Im}

\newcommand{\Ext}[3]{\mbox{Ext}^1_{#1}\,(#2,#3)}

\newcommand{\Exti}[4]{\mbox{Ext}^{#1}_{#2}\,(#3,#4)}

\newcommand{\Tor}[3]{\mbox{Tor}_1^{#1}\,(#2,#3)}

\DeclareMathOperator{\HomOp}{Hom}

\newcommand{\Hom}[3]{\HomOp_{#1}(#2,#3)}

\newcommand{\Dcal}{\ensuremath{\mathcal{D}}}

\newcommand{\Ucal}{\ensuremath{\mathcal{U}}}

\newcommand{\Qcal}{\ensuremath{\mathcal{Q}}}

\newcommand{\Ycal}{\ensuremath{\mathcal{Y}}}
\newcommand{\Xcal}{\ensuremath{\mathcal{X}}}

\newcommand{\p}{\ensuremath{\mathbf{p}}}

\newcommand{\tube}{\ensuremath{\mathbf{t}}}

\newcommand{\Tria}{\mathrm{Tria\,}}
\theoremstyle{plain}
\newtheorem{thm}{Theorem}[section]
\newtheorem{prop}[thm]{Proposition}
\newtheorem{lem}[thm]{Lemma}
\newtheorem{lemma}[thm]{Lemma}
\newtheorem{cor}[thm]{Corollary}
\newtheorem{ex}[thm]{Example}
\theoremstyle{definition}

\theoremstyle{remark}
\newtheorem*{rem}{Remark}

\newcommand{\fatbox}{}

\newcommand{\ra}{\rightarrow}
\newcommand{\lra}{\longrightarrow}
\newcommand{\les}{\leqslant}
\newcommand{\ges}{\geqslant}
\newcommand{\ten}{\otimes}
\newcommand{\lten}{\overset{\boldmath{L}}{\ten}}

\newcommand{\mcy}{\mathcal{Y}}
\newcommand{\mcx}{\mathcal{X}}

\begin{document}
\begin{center}
 {\bf Recollements and tilting objects}
\bigskip

{\sc Lidia Angeleri  H\" ugel, Steffen Koenig, Qunhua Liu}

\bigskip

(Version of \today)
 \end{center}

\address{Lidia Angeleri  H\" ugel\\ Dipartimento di Informatica -
Settore Matematica
\\ Universit\`a degli Studi di Verona
\\ Strada Le Grazie 15 - Ca' Vignal 2
\\ I - 37134 Verona, Italy}
\email{ lidia.angeleri@univr.it}
\address{Steffen Koenig, Qunhua Liu\\ Mathematisches Institut der
Universit\"at zu K\"oln \\ Weyertal 86-90 \\ 50931 K\"oln, Germany}
\email{skoenig@math.uni-koeln.de, qliu@math.uni-koeln.de}

\date{\today}


\bigskip

\begin{quote}
{\footnotesize {\bf Abstract.} We study connections between recollements of the derived category D(Mod$R$) of a ring $R$ and
tilting theory.
We first  provide constructions of
tilting objects from given recollements, recovering several different
results from the literature. Secondly, we  show how to construct a
recollement from a  tilting module of projective
dimension one. By \cite{NS1}, every recollement of D(Mod$R$) is associated to a differential graded homological epimorphism $\lambda:R\to S$. We will focus on the case where $\lambda$ is  a homological ring epimorphism or even a universal localization.
Our results will be employed in a forthcoming paper in order to investigate stratifications of D(Mod$R$).}

\end{quote}

\vspace{4ex}

\section*{Introduction}

Recollements of triangulated categories are 'exact sequences' of triangulated
categories, which describe the middle term by a triangulated subcategory and
a triangulated quotient category. Recollements have  first been defined by
Beilinson, Bernstein and Deligne \cite{BBD} in a geometric context, where
stratifications of spaces imply recollements of derived categories of sheaves,
by using derived versions of
Grothendieck's six functors (which conveniently get axiomatized by
the concept of recollement). As certain derived categories of perverse sheaves
are equivalent to derived categories of modules over blocks of the
Bernstein-Gelfand-Gelfand category $\mathcal O$, recollements do exist for
the corresponding algebras as well. Here, the stratification provided by
iterated recollements, is by derived categories of vector spaces. This is
one of the fundamental, and motivating,
properties of quasi-hereditary algebras, introduced
by Cline, Parshall and Scott (see \cite{PS}).

\medskip

The first examples of recollements of derived categories of rings
have been produced by direct constructions, using derived functors
of known functors on abelian level. Subsequently, a necessary and
sufficient criterion has been given \cite{K} for a (bounded)
derived module category of an algebra to admit a recollement, with
subcategory and quotient category again being derived module
categories of rings. This criterion is formulated in terms of two
exceptional objects that fully describe the recollement. Later on,
the criterion has been extended and modified so as to cover
derived categories of differential graded algebras and unbounded
derived categories as well and to work for any differential graded
ring \cite{J,NS1}. All these results characterize the existence of a
recollement in terms of two exceptional objects. In the special
case of the quotient or the subcategory being zero, one
exceptional object is zero and the other is a tilting complex,
that is, one recovers Morita theory of derived categories. While
in this special case, the role of the tilting complex is very
natural in the context of tilting theory, little is known about
connections between recollements of derived module categories and
tilting theory. The aim of this article is to start exploring such
potential connections. We  will first provide constructions of
tilting objects from given recollements. Our constructions will be
general enough to cover quite a few, and rather diverse,
situations studied in the literature (usually without mentioning
recollements). Conversely, we will show how to construct a
recollement from a classical tilting module (of projective
dimension one); in this way we will extend results in \cite{AHS}
and put them into a general framework.

\medskip

In the first section we will collect existence results and
categorical methods to construct recollements. The second section
leads to the first main result, Theorem \ref{constructtogen} and
its variation Theorem \ref{constructto} (for a situation
satisfying some finiteness conditions), which construct a tilting
object from the two exceptional objects describing a recollement;
the axioms of a recollement imply that there are no morphisms
between the two exceptional objects in one direction, and we also
assume that morphism in the opposite direction are concentrated in
at most two degrees. The subsequent section three applies the
first main result in quite diverse situations, thus recovering and
re-interpreting various results from the literature. In the fourth
section we start with a classical or a large tilting module of
projective dimension one over any ring, and construct a
recollement from it. The main result, Theorem \ref{recoll},
describes both the subcategory and the quotient category in such a
recollement. The latter is a derived module category in the
classical case; the former is shown to be equivalent to a derived
module category if and only if a certain universal localization is
a homological epimorphism. Examples of such situations are given
in the final section; some of these examples also illustrate
differences between various technical terms used in developing the
theory. In an appendix, we provide a construction for reflections
in triangulated categories.

\medskip

In the subsequent article \cite{AKL}, we will be strongly using the results of
the present article to address a basic and so far completely open question
about recollements: Is there a Jordan-H\"older theorem for derived categories?
In other words,
is there an existence and uniqueness result for iterated recollements
(that is, for stratifications of derived categories)?
We will show by various examples of 'exotic stratifications' that the answer
(and the validity of such a Jordan H\"older theorem)
depends very much on the choice of triangulated categories (such as
derived categories of algebras or of differential graded algebras or other
triangulated categories). Moreover, we will provide positive answers;
in particular, we will
prove a Jordan-H\"older theorem for bounded derived categories of artinian
hereditary rings and thus also for all piecewise hereditary algebras. Here,
crucial use will be made in particular of Theorem \ref{recoll}, which will
allow to identify the end terms of certain recollement situations as derived
module categories.
We will also  discuss when hereditary rings are derived simple.

\bigskip

{\bf Acknowledgements:} We thank Nan Gao for simplifying our proof
of Lemma \ref{perpendicular}.

First named author acknowledges partial support from  MIUR, PRIN-2008 "Anelli, algebre, moduli e categorie", and from Progetto di Ateneo CPDA071244 of the University of Padova.

\bigskip

\section{Recollements and localizations}

In this section, recollements are defined and various criteria for the
existence of recollements are discussed.

Throughout this paper, $\Dcal$ denotes a triangulated category with small
coproducts (that is, coproducts indexed over a set), and $[1]$ denotes the
shift functor.


\medskip
\subsection{Recollements.}
Let $\Xcal, \Ycal$ be triangulated categories. $\Dcal$ is said to be  a
{\em recollement} of $\Xcal$ and $\Ycal$ if there are six triangle functors
as in the following diagram
\vspace{0.3cm}
\begin{center}
\unitlength0.3cm
\begin{picture}(15,2)
\put(-0.5,0){{\Ycal}}
\put(6,0){{\Dcal}}
\put(12.9,0){{\Xcal}}
\put(3.2,-0.3){\oval(4,2.1)[b]}
\put(3.2,0.4){\oval(4,2.9)[t]}
\put(1.2,-0.3){\vector(0,1){0.3}}
\put(1.2,0.8){\vector(0,-1){0.3}}
\put(1.6,0.2){\vector(1,0){3.3}}
\put(2.7,2.3){$\scriptstyle i^\ast$}
\put(2.7,-1.33){$\scriptstyle i^!$}
\put(2,0.5){$\scriptstyle i_\ast=i_!$}
\put(10.1,-0.3){\oval(4,2.1)[b]}
\put(10.1,0.4){\oval(4,2.9)[t]}
\put(8.1,-0.3){\vector(0,1){0.3}}
\put(8.1,0.8){\vector(0,-1){0.3}}
\put(8.4,0.2){\vector(1,0){3.3}}
\put(9.5,2.3){$\scriptstyle j_\ast$}
\put(9.5,-1.15){$\scriptstyle j_!$}
\put(9,0.5){$\scriptstyle j^!=j^\ast$}
\end{picture}
\end{center}
such that
\begin{enumerate}
\item
$(i^\ast,i_\ast)$,\,$(i_!,i^!)$,\,$(j_!,j^!)$ ,\,$(j^\ast,j_\ast)$  are adjoint
pairs;

\smallskip

 \item
 $i_\ast,\,j_\ast,\,j_!$  are full embeddings;

 \smallskip

 \item  $i^!\circ j_\ast=0$ (and thus also $j^!\circ i_!=0$ and
$i^\ast\circ j_!=0$);

 \smallskip

 \item
 for each $C\in \Dcal$ there are triangles $$i_! i^!(C)\to
C\to j_\ast j^\ast (C)\to $$
 $$j_! j^! (C)\to C\to i_\ast i^\ast(C)\to$$
\end{enumerate}

\bigskip

Recollements are closely related to localization, which will be discussed
below.

\medskip

\subsection{Bousfield localization.}\label{localiz}
A triangle functor $L:\Dcal\to \Dcal$ is said to be a
{\em localization functor}
if there is a natural transformation $\eta:{\mathrm Id}\to L$ such that for all
$X\in\Dcal$\\
(i) $L \circ \eta_{X}=\eta_{L(X)}$ , and\\
(ii) $\eta$ induces an isomorphism $L(X)\cong L^2(X)$.

\medskip

Such a localization functor determines a full subcategory $\Xcal$ of $\Dcal$
whose objects are precisely the $X\in\Dcal$ such that $L(X)=0$.
Subcategories of
$\Dcal$ arising in this way are called {\em localizing subcategories}.

\medskip

Note that $\Xcal$ is a thick subcategory of $\Dcal$, so we can form the
quotient
category $\Dcal/\Xcal$, see \cite{Ve}.
We consider the quotient functor  $$\pi: \Dcal\to\Dcal/\Xcal$$ and we denote
by
$\Ycal$  the  right orthogonal class of $\Xcal$ given by all objects
$Y\in\Dcal$
such that $\Hom {\Dcal}{X}{Y} =0$ for all $X\in\Xcal$.

The following statements hold true
(see  e.g. \cite[1.6]{AJS1})\label{localizing}
\begin{enumerate}
\item
The  functor $\pi:\Dcal\to \Dcal/\Xcal$  induces an equivalence
$\pi\circ {{\rm inc}_{\Ycal}}:\Ycal\to \Dcal/\Xcal$ with inverse
$\rho$.
\item
The functor ${{\rm inc}_{\Ycal}}$ has a left adjoint $q=\rho\circ \pi$.
\item The functor ${{\rm inc}_{\Xcal}}$ has a right adjoint $a$.
\end{enumerate}

\medskip

We thus obtain triangle functors as in the following diagram:

\vspace{0.3cm}
\begin{center}
\unitlength0.3cm
\begin{picture}(15,2)
\put(-0.5,0){{\Ycal}}
\put(6,0){{\Dcal}}
\put(12.6,0){{\Xcal}}
\put(3.2,0.4){\oval(4,2.9)[t]}
\put(1.2,0.8){\vector(0,-1){0.3}}
\put(1.6,0.2){\vector(1,0){3.3}}
\put(2.7,2.3){$\scriptstyle q$}
\put(2.7,0.5){$\scriptstyle {\rm inc}$}
\put(10,-0.3){\oval(4,2.1)[b]}
\put(8,-0.3){\vector(0,1){0.3}}
\put(8.4,0.2){\vector(1,0){3.3}}
\put(9.5,-1.15){$\scriptstyle {\rm inc}$}
\put(9.5,0.5){$\scriptstyle a$}
\end{picture}
\end{center}

\smallskip

\smallskip

where $q\circ {\rm inc}_\Xcal=0$ and $a\circ {\rm inc}_\Ycal=0$.

\bigskip

Note that the localization functor $L$ preserves small coproducts
if and only if the category $\Ycal$ is closed under small
coproducts. In this case the localizing subcategory $\Xcal$ is
said to be a {\em smashing subcategory}, and \label{smashing}
there even is a {recollement} \vspace{0.3cm}
\begin{center}
\unitlength0.3cm
\begin{picture}(15,2)
\put(-0.5,0){{\Ycal}}
\put(6,0){{\Dcal}}
\put(13,0){{\Xcal}}
\put(3.2,-0.3){\oval(4,2.1)[b]}
\put(3.2,0.4){\oval(4,2.9)[t]}
\put(1.2,-0.3){\vector(0,1){0.3}}
\put(1.2,0.8){\vector(0,-1){0.3}}
\put(1.6,0.2){\vector(1,0){3.3}}
\put(2.7,2.3){$\scriptstyle q$}
\put(2.7,-1.35){$\scriptstyle b$}
\put(2.7,0.5){$\scriptstyle {\rm inc}$}
\put(10.4,-0.3){\oval(4,2.1)[b]}
\put(10.4,0.4){\oval(4,2.9)[t]}
\put(8.4,-0.3){\vector(0,1){0.3}}
\put(8.4,0.8){\vector(0,-1){0.3}}
\put(8.8,0.2){\vector(1,0){3.3}}
\put(9.9,2.3){$\scriptstyle j$}
\put(9.9,-1.15){$\scriptstyle {\rm inc}$}
\put(9.9,0.5){$\scriptstyle a$}
\end{picture}
\end{center}
More precisely,
\begin{enumerate}
\item
the functor ${\rm inc}_{\Ycal}$ has a right adjoint $b$,

\smallskip

\item the functor $a$ has a right  adjoint $j$,

\smallskip

 \item
 $j$ is a full embedding, and  $b\circ j=0$;

 \smallskip

 \item
 for each $C\in \Dcal$ there are triangles $${\rm inc}_{\Ycal} b (C)\to C\to j a
(C)\to $$
 $${\rm inc}_{\Xcal} a (C)\to C\to {\rm inc}_{\Ycal} q (C)\to$$
\end{enumerate}

\medskip

For details on the correspondence between smashing subcategories and recollements
we refer to \cite[4.4.14, 4.2.4, 4.2.5]{N}, \cite{NS1}.

\bigskip

Let us now turn to our main example.

\medskip

\subsection{The derived category of a ring.} Let $R$ be a ring, and let \M\
be the category of all right
$R$-modules.
We denote by $\Dcal(R)$ the unbounded derived category of $\M$. The category $\M$
is identified with the subcategory of $\Dcal(R)$ consisting of the stalk
complexes concentrated in degree zero. Of course,  every module $M$ is
quasi-isomorphic
to the complex given by a projective resolution of $M$.

\medskip

\subsection{Generators, compact objects, tilting objects.} \label{generators}
Given a class of objects $\Qcal$ in $\Dcal$,
the smallest full triangulated subcategory of $\Dcal$ which contains $\Qcal$ and
is closed under small coproducts is denoted by
$\Tria\Qcal$ (note that some
authors use the notation $\Tria^{{\scriptscriptstyle\coprod}}\Qcal$).
If $\Qcal$ consists just of one object $Q$, we write
$\Tria {Q}$.

\smallskip

The triangulated category $\Dcal$ satisfies the \emph{principle of
infinite d\'evissage} (with respect to $\Qcal$) if
$\Dcal=\Tria\Qcal$. In this case, $\Dcal$ is \emph{generated} by
$\Qcal$, that is,  an object $M$ of $\Dcal$ is zero whenever $\Hom
{\Dcal}{Q[n]}{M}=0$ for every object $Q$ of $\Qcal$ and every
$n\in\Z$. Sometimes also the converse holds true. For example, if
$\Ycal$ is a full triangulated subcategory of $\Dcal$ generated by
$\Qcal$ and $\Tria\Qcal$ is an aisle in $\Dcal$ contained in
$\Ycal$, then $\Ycal=\Tria\Qcal$, see \cite[4.3.5 and 4.3.6]{N}, \cite{NS1}.

\smallskip

An object $P$ of $\Dcal$ is said to be \emph{compact} if the functor
$\Hom {\Dcal}{{P}}{-}$ preserves small coproducts.
Furthermore,
$P$ is said to be \emph{self-compact} if the restricted functor
$\Hom {\Dcal}{{P}}{-}\mid_{\Tria{P}}$ preserves small coproducts.

\smallskip

It is well known that a complex ${P^\cdot}\in\Dcal(R)$ is compact
if and only if it is quasi-isomorphic to a bounded complex
consisting of finitely generated projective modules. In
particular,  the compact objects of $\M$ are precisely the modules
in \m\ of finite projective dimension. Here  \m\ denotes the
subcategory of $\M$ given by all modules possessing a projective
resolution consisting of finitely generated modules.

\medskip

 An object
$T$ in $\Dcal$ is called  {\em exceptional} (or a partial tilting object) if $T$
has no
self extensions, i.e. $\text{Hom}_\Dcal(T,T[k])=0$ for all nonzero integers
$k$.
Furthermore,
$T$ is called a {\em tilting} object if it is compact, exceptional, and $\Dcal$
is  generated by $T$.

\medskip

We will frequently use the following result due to Keller.

\smallskip

{\bf Theorem.} \cite[Theorem 8.5]{KHandbook}, \cite{Ke} {\it Let
$R$ be a ring, and let $\Dcal$ be a full triangulated subcategory
of $\Dcal(R)$ closed under coproducts. If $T$ is a compact
generator of $\Dcal$, then  there is a differential graded algebra
$E=\mathbf{R}\Hom{}{T}{T}$ with homology
 $H^\ast(E)\cong \bigoplus_{i\in\Z}\Hom{\Dcal}{T}{T[i]}$ such that the functor
$-\otimes^{\mathbf L}_E T: \Dcal(E)\to \Dcal$ is a  triangle equivalence.
}

\medskip

\subsection{Localizing subcategories generated by a set.}\label{bousfield}
By results of Bousfield and Neeman, every set $\Qcal$ of compact
objects in $\Dcal(R)$ defines a smashing subcategory $\Tria\Qcal$
and therefore a recollement of $\Dcal(R)$ (see e.g. \cite[4.4.16
and 4.4.3]{N}). We will   often work under weaker assumptions and will need a result from \cite{AJS1}
stating that {\em any   set} of objects   in $\Dcal(R)$ gives rise
to a localizing subcategory.

\medskip

{\bf Theorem \cite[4.5]{AJS1}} {\it Let $\Qcal$ be a set of
objects in $\Dcal(R)$. Set $\Xcal=\Tria\Qcal$, and let
$\Ycal=\Ker\Hom {\Dcal}{\Xcal}{-}$ be the right orthogonal class.
Then $\Xcal$ is a localizing subcategory of $\Dcal(R)$,
and $\Ycal$ consists of the objects $Y^\cdot\in\Dcal(R)$ such that
$\Hom {\Dcal(R)}{Q^\cdot[n]}{Y^\cdot} =0$ for all $Q^\cdot\in\Qcal$ and $n\in\Z$.
If $\Qcal$ consists of compact objects, then   $\Xcal$ is even a smashing
subcategory.}

\bigskip

\subsection{Recollements induced by single objects} \label{singleobjects}
The following result was first proved by the second named author
for bounded derived categories \cite{K}, and it was then  further
developed by several authors \cite{J,N} (note that in \cite{K} a
condition has been misstated, see \cite{NS2} for a discussion). The
versions of this result in \cite{K} and in \cite{J} are assuming
that all triangulated categories are derived categories of
(differential graded) rings; therefore, the exceptional objects
that appear there are images of two of the rings. The exceptional
objects appearing in the following version are, in general
different, even if all categories are derived categories of rings.

\medskip

{\bf Theorem. (\cite[5.2.9]{N}, \cite{NS1})}\label{single} {\it
The derived category $\Dcal(R)$ of a ring $R$ is a recollement of
derived categories of rings if and only if
there are  objects $T_1,T_2\in\Dcal(R)$ such that\\
(i) $T_1$ is compact and exceptional,\\
(ii) $T_2$ is self-compact and exceptional,\\
(iii) $\Hom{\Dcal}{T_1[n]}{T_2}=0$ for all $n\in\Z$,\\
(iv) $\{T_1,T_2\}$ generates $\Dcal(R)$.}

\medskip

We will need the following ``non-compact version'' of this theorem.

\medskip

{\bf Theorem.}\label{singlenoncompact} {\it Assume that $\Dcal$ has a compact
generator $R$. Then the following statements are equivalent.
\begin{enumerate}
\item $\Dcal$ is a recollement of triangulated categories generated by a single object.
\item There is an object $T_1\in\Dcal$ such that $\Tria T_1$ is a smashing
subcategory of $\Dcal$.
\item There is an object $T_1\in\Dcal$ such that $\Ker\Hom{\Dcal}{\Tria T_1}{-}$ is closed under coproducts.
\item There are objects $T_1,T_2\in\Dcal$ such that\\
(i) $\Ker\Hom{\Dcal}{\Tria T_1}{-}$ is closed under coproducts,\\
(ii) $T_2$ is self-compact,\\
(iii) $\Hom{\Dcal}{T_1[n]}{T_2}=0$ for all $n\in\Z$,\\
(iv) $\{T_1,T_2\}$ generates $\Dcal$.
\end{enumerate}
}

\begin{proof} (1)$\Rightarrow$(2): Condition (1) implies the existence of a
smashing subcategory $\Xcal$ generated by an object $T_1$. We have
just seen in \ref{bousfield} that $\Tria T_1$ is a localizing
subcategory of $\Dcal$ (i.e. an aisle in $\Dcal$) which is contained in
$\Xcal$. So, we infer from \cite[4.3.6]{N} that $\Xcal=\Tria T_1$.

By \ref{localiz} and  \ref{bousfield}, the
conditions (2) and (3) are equivalent.
(4)$\Rightarrow$(3) is clear.

It remains to show (3)$\Rightarrow$(4),(1): It follows from
condition (3) that there is a recollement \vspace{0.15cm}
\begin{center}
\unitlength0.3cm
\begin{picture}(15,2)
\put(-0.5,0){{\Ycal}}
\put(6,0){{\Dcal}}
\put(13,0){{$\Xcal=\Tria T_1$}}
\put(3.2,-0.3){\oval(4,2.1)[b]}
\put(3.2,0.4){\oval(4,2.9)[t]}
\put(1.2,-0.3){\vector(0,1){0.3}}
\put(1.2,0.8){\vector(0,-1){0.3}}
\put(1.6,0.2){\vector(1,0){3.3}}
\put(2.7,2.3){$\scriptstyle q$}
\put(2.7,-1.35){$\scriptstyle b$}
\put(2.7,0.5){$\scriptstyle {\rm inc}$}
\put(10.4,-0.3){\oval(4,2.1)[b]}
\put(10.4,0.4){\oval(4,2.9)[t]}
\put(8.4,-0.3){\vector(0,1){0.3}}
\put(8.4,0.8){\vector(0,-1){0.3}}
\put(8.8,0.2){\vector(1,0){3.3}}
\put(9.9,2.3){$\scriptstyle j$}
\put(9.9,-1.15){$\scriptstyle {\rm inc}$}
\put(9.9,0.5){$\scriptstyle a$}
\end{picture}
\end{center}
\medskip
 \noindent and by \cite[4.3.6, 4.4.8]{N}, \cite{NS1}, the
compact generator $R$ of $\Dcal$ is mapped by $q$ to a compact
generator $T_2=q(R)$ of $\Ycal$. As above, we infer $\Ycal=\Tria
T_2$, and we immediately verify (ii) and (iii). Finally, condition
(iv) follows from the existence of  triangles ${\rm inc}_{\Xcal} a
(C)\to C\to {\rm inc}_{\Ycal} q (C)\to $ with $q(C)\in\Tria T_2$
and $a(C)\in\Tria T_1$ for each object $C\in\Dcal$.
\end{proof}

In the case when $\Dcal=\Dcal(R)$ and  $T_1$ is compact and exceptional, we provide a
construction of the object $T_2=q(R)$ in the Appendix. More
precisely, we construct the  $\Ycal$-reflection  $M\ra
q(M)$ of $M$ for those $M\in\Dcal$ such that
$\Hom{\Dcal}{T_1}{M[i]}=0$ for sufficiently large $i$.

\bigskip

 Here is another source of examples for recollements.

\subsection{Homological ring epimorphisms.}\label{ringepi} Let $\lambda\colon
R\rightarrow S$ be a ring epimorphism, that is, an epimorphism in the category of
rings.
 Following  Geigle and Lenzing \cite{GL}, we say that
   $\lambda$ is  a  \emph{homological ring
epimorphism} if $\mbox{Tor} _i^R(S,S)=0$ for all $i>0$. Note that this holds true
if and only if the restriction functor $\lambda_\ast:\Dcal(S)\to\Dcal(R)$
induced  by $\lambda$ is fully faithful \cite[4.4]{GL},\cite[5.3.1]{N}. As shown
in \cite[Section 5.3]{N},\cite{NS1}, we then obtain  a recollement
\vspace{0.3cm}
\begin{center}
\unitlength0.3cm
\begin{picture}(7,2)
\put(-1.8,0){{$\Dcal(S)$}}
\put(5.6,0){{$\Dcal(R)$}}
\put(13.3,0){{$\Tria X$}}
\put(3.2,-0.3){\oval(4,2.1)[b]}
\put(3.2,0.4){\oval(4,2.9)[t]}
\put(1.2,-0.3){\vector(0,1){0.3}}
\put(1.2,0.8){\vector(0,-1){0.3}}
\put(1.6,0.2){\vector(1,0){3.3}}
\put(2.7,2.1){$\scriptstyle F$}
\put(2.7,-1.28){$\scriptstyle G$}
\put(2.7,0.5){$\scriptstyle \lambda_\ast$}
\put(10.7,-0.3){\oval(4,2.1)[b]}
\put(10.7,0.4){\oval(4,2.9)[t]}
\put(8.7,-0.3){\vector(0,1){0.3}}
\put(8.7,0.8){\vector(0,-1){0.3}}
\put(9.1,0.2){\vector(1,0){3.3}}
\put(10.2,2.3){$\scriptstyle $}
\put(10.2,-1.15){$\scriptstyle $}
\put(10.2,0.5){$\scriptstyle \tau$}
\end{picture}
\end{center}

\bigskip

\noindent
where $F=-\otimes_R^{\mathbf L} S$ is the derived tensor product,
$G={\mathbf R} \Hom R{S}{-}$ is the derived Hom-functor, $X$ is the object
occurring in the triangle $$X\to R\stackrel{\lambda}{\to} S\to$$
and $\tau=-\otimes_R^{\mathbf L} X$.
This also follows
from   \cite[Theorem 2.4 (1)]{PS} (which proves that a `partial'
recollement can be completed).

\medskip

There is also a converse result: by  \cite[5.4.4]{N}, \cite{NS1}, every recollement of
$\Dcal(R)$ is associated to a differential graded homological epimorphism
$\lambda:R\to S$.
    In this paper, we will focus on the case of $\lambda$ being a homological
ring epimorphism.

    Following \cite{GP}, we will say that two ring epimorphisms
$\lambda :R\to S$ and $\lambda' :R\to S'$ are equivalent if there is
a ring isomorphism $\psi :S\to S'$ such that $\lambda'=\psi
\lambda$. The   equivalence classes
with respect to this equivalence relation are called \emph{epiclasses}.

Moreover, we will say that two recollements

\vspace{0.3cm}
\begin{center}
\unitlength0.3cm
\begin{picture}(30,2)
\put(-0.5,0){{\Ycal}}
\put(6,0){{\Dcal}}
\put(12.9,0){{\Xcal}}
\put(3.2,-0.3){\oval(4,2.1)[b]}
\put(3.2,0.4){\oval(4,2.9)[t]}
\put(1.2,-0.3){\vector(0,1){0.3}}
\put(1.2,0.8){\vector(0,-1){0.3}}
\put(1.6,0.2){\vector(1,0){3.3}}
\put(2.7,2.3){$\scriptstyle $}
\put(2.7,-1.33){$\scriptstyle$}
\put(2.5,0.5){$\scriptstyle i_\ast$}
\put(10.1,-0.3){\oval(4,2.1)[b]}
\put(10.1,0.4){\oval(4,2.9)[t]}
\put(8.1,-0.3){\vector(0,1){0.3}}
\put(8.1,0.8){\vector(0,-1){0.3}}
\put(8.4,0.2){\vector(1,0){3.3}}
\put(9.5,2.3){$\scriptstyle j_\ast$}
\put(9.5,-1.15){$\scriptstyle j_!$}
\put(9,0.5){$\scriptstyle $}

\put(15,0){\text{and}}

\put(18.5,0){$\Ycal'$} \put(25,0){{\Dcal}} \put(31.9,0){$\Xcal'$}
\put(22.2,-0.3){\oval(4,2.1)[b]} \put(22.2,0.4){\oval(4,2.9)[t]}
\put(20.2,-0.3){\vector(0,1){0.3}}
\put(20.2,0.8){\vector(0,-1){0.3}}
\put(20.6,0.2){\vector(1,0){3.3}}
\put(21.7,2.3){$\scriptstyle $}
\put(21.7,-1.33){$\scriptstyle$}
\put(21.5,0.5){$\scriptstyle i'_\ast$}
\put(29.1,-0.3){\oval(4,2.1)[b]}
\put(29.1,0.4){\oval(4,2.9)[t]}
\put(27.1,-0.3){\vector(0,1){0.3}}
\put(27.1,0.8){\vector(0,-1){0.3}}
\put(27.4,0.2){\vector(1,0){3.3}}
\put(28.5,2.3){$\scriptstyle j'_\ast$}
\put(28.5,-1.15){$\scriptstyle j'_!$}
\put(28,0.5){$\scriptstyle $}
\end{picture}
\end{center}

\bigskip

\noindent are \emph{equivalent} if the essential images of
$i_\ast$ and $i'_\ast$, of $j_\ast$ and $j'_\ast$, and of $j_!$
and $j'_!$ coincide, respectively.

    The following observation is implicit in \cite{N, NS1}.

    \bigskip

\noindent
{\bf Proposition.}
{\it Let $R$ be a ring and $\Dcal = D(R)$ its derived category.
Then there is a bijection between the epiclasses of homological ring
epimorphisms starting in $R$ and the equivalence classes of those recollements
\vspace{0.3cm}
\begin{center}
\unitlength0.3cm
\begin{picture}(13,2)
\put(-0.5,0){{\Ycal}}
\put(6,0){{\Dcal}}
\put(12.9,0){{\Xcal}}
\put(3.2,-0.3){\oval(4,2.1)[b]}
\put(3.2,0.4){\oval(4,2.9)[t]}
\put(1.2,-0.3){\vector(0,1){0.3}}
\put(1.2,0.8){\vector(0,-1){0.3}}
\put(1.6,0.2){\vector(1,0){3.3}}
\put(2.7,2.3){$\scriptstyle i^\ast$}
\put(2.7,-1.33){$\scriptstyle $}
\put(2.5,0.5){$\scriptstyle i_\ast$}
\put(10.1,-0.3){\oval(4,2.1)[b]}
\put(10.1,0.4){\oval(4,2.9)[t]}
\put(8.1,-0.3){\vector(0,1){0.3}}
\put(8.1,0.8){\vector(0,-1){0.3}}
\put(8.4,0.2){\vector(1,0){3.3}}
\put(9.5,2.3){$\scriptstyle $}
\put(9.5,-1.15){$\scriptstyle $}
\put(9,0.5){$\scriptstyle $}
\end{picture}
\end{center}

\bigskip

\noindent
for which $i^\ast(R)$ is an exceptional object of $\Ycal$.}

\bigskip

\begin{proof}
Let $\lambda:R\to S$ be a  homological ring epimorphism, and consider the
recollement induced by $\lambda$ as above. Then the image of $R$ under the
functor $F=-\otimes_R^{\mathbf L} S$ is isomorphic to $S$ and  thus  an
exceptional object of $\Dcal(S)$.

Conversely, take  a recollement as in the Proposition, for which $Q=i^\ast(R)$
is an exceptional object of $\Ycal$. Note that $Q$ is a compact generator of
$\Ycal$ by \cite[4.3.6 and 4.4.8]{N}, whence a compact tilting object in
$\Ycal$. By Keller's theorem in \ref{generators} there is a
differential graded algebra $E=\mathbf{R}\Hom{}{Q}{Q}$ having homology
concentrated in zero and $H^0(E)\cong \Hom{\Ycal}{Q}{Q}\cong
\text{End}_{\Dcal(R)}\,i_\ast(Q)$, such that the functor
$-\otimes^{\mathbf L}_E Q$ defines a  triangle equivalence between  the derived
category of $E$ and $\Ycal$.  Setting $S=\text{End}_{\Dcal(R)}\,i_\ast(Q)$ we
obtain a full embedding $\iota:\Dcal(S)\to \Dcal(R)$.

Note that $H^n(i_\ast (Q))\cong \Hom{\Dcal(R)}{R}{i_\ast(Q)[n]} \cong
\Hom{\Ycal}{Q}{Q[n]}$,
hence  $i_\ast (Q)$ has homology concentrated in zero, and
$H^0(i_\ast (Q))\cong S$.

Moreover, the unit $\eta$ of the adjoint pair $(i^\ast,i_\ast)$ yields a
$\Ycal$-reflection $$\eta_R:R\to i_\ast i^\ast(R)=i_\ast(Q),$$ that is,
$\Hom {\Dcal(R)}{\eta_R}{i_\ast(Y)}:\Hom {\Dcal(R)}{i_\ast(Q)}{i_\ast(Y)}\to
\Hom {\Dcal(R)}{R}{i_\ast(Y)}$ is a bijection for every $Y\in\Ycal$.
This allows to define a ring homomorphism $\lambda:R\to S$ by associating to
any element $r\in R$ the left multiplication $m_r: R\to R, x\mapsto rx$ and
setting
$\lambda(r)= i_\ast i^\ast(m_r):i_\ast(Q)\to i_\ast(Q)$.

In this way, $S$ becomes a right $R$-module, that is, a complex concentrated
in zero, which is quasi-isomorphic to  $i_\ast(Q)$.  It follows that the
restriction functor $\lambda_\ast:\Dcal(S)\to\Dcal(R)$ induced  by $\lambda$
coincides with the full embedding $\iota:\Dcal(S){\to}\Dcal(R)$, showing that
$\lambda$ is a homological ring epimorphism.

Now it is clear how to define the stated bijective correspondence.  \fatbox
\end{proof}

    \bigskip

    \subsection{Universal localization.}\label{universal}
    Finally, we focus on a special kind of homological ring epimorphisms.

\smallskip

{\bf Theorem.} {\cite[Theorem~4.1]{Schofieldbook}}
\label{def:universallocalization} {\it Let  $\Sigma$ be a set of
morphisms between finitely generated projective right $R$-modules.
Then there are a ring $R_\Sigma$ and a morphism of rings
$\lambda\colon R\rightarrow R_\Sigma$ such that
\begin{enumerate}
\item $\lambda$ is \emph{$\Sigma$-inverting,} i.e. if
$\alpha\colon P\rightarrow Q$ belongs to  $\Sigma$, then
$\alpha\otimes_R 1_{R_\Sigma}\colon P\otimes_R R_\Sigma\rightarrow
Q\otimes_R R_\Sigma$ is an isomorphism of right
$R_\Sigma$-modules, and
\item $\lambda$ is \emph{universal
$\Sigma$-inverting}, i.e. if $S$ is a ring such that there exists
a $\Sigma$-inverting morphism $\psi\colon R\rightarrow S$, then
there exists a unique morphism of rings $\bar{\psi}\colon
R_\Sigma\rightarrow S$ such that $\bar{\psi}\lambda=\psi$.
\end{enumerate}
}

\medskip

The morphism $\lambda\colon R\rightarrow R_\Sigma$ is a ring epimorphism
with  $\Tor R{R_\Sigma}{R_\Sigma}=0.$ It
 is called the \emph{universal localization of $R$ at
$\Sigma$}.

\medskip

Let now $\mathcal{U}$ be a set of  finitely presented right $R$-modules of
projective dimension one.
For each $U\in\mathcal{U},$ consider a morphism $\alpha_U$ between
finitely generated projective right $R$-modules such that
$$0\to P\stackrel{{\alpha_U}}{\to} Q\to U\to 0$$
 We will denote by
$R_{\mathcal{U}}$ the universal localization of $R$ at
$\Sigma=\{\alpha_U\mid U\in\mathcal{U}\}.$ In fact,   $R_{\mathcal
U}$ does not depend on the class $\Sigma$  chosen,
cf.~\cite[Theorem~0.6.2]{Cohnfreerings}, and we will also call it
the \emph{universal localization of $R$ at ${\mathcal{U}}$}.

\medskip

In general, a universal localization  need not be a homological ring
epimorphism, see \cite{NRS} and Example \ref{nothom}.
Universal localizations with this stronger homological property were studied
by Neeman and Ranicki. We will need the following result, which is a
combination of some of their results in \cite{NR}.

\smallskip

{\bf Theorem.} \label{univloc}
{\it Let  $\mathcal{U}$ be a set of  finitely presented right $R$-modules of
projective dimension one. Assume that the universal localization
$\lambda:R\to R_{\Ucal}$ is a homological ring epimorphism. Then there is a
recollement
\vspace{0.3cm}
\begin{center}
\unitlength0.3cm
\begin{picture}(15,2)
\put(-2.8,0){{$\Dcal(R_{\Ucal})$}}
\put(5.6,0){{$\Dcal(R)$}}
\put(13.3,0){{$\Tria\Ucal$}}
\put(3.2,-0.3){\oval(4,2.1)[b]}
\put(3.2,0.4){\oval(4,2.9)[t]}
\put(1.2,-0.3){\vector(0,1){0.3}}
\put(1.2,0.8){\vector(0,-1){0.3}}
\put(1.4,0.2){\vector(1,0){3.5}}
\put(2.7,2.2){$\scriptstyle F$}
\put(2.7,-1.2){$\scriptstyle G$}
\put(2.7,0.5){$\scriptstyle {\rm \lambda_\ast}$}
\put(10.7,-0.3){\oval(4,2.1)[b]}
\put(10.7,0.4){\oval(4,2.9)[t]}
\put(8.7,-0.3){\vector(0,1){0.3}}
\put(8.7,0.8){\vector(0,-1){0.3}}
\put(8.9,0.2){\vector(1,0){3.5}}
\put(10.2,2.3){$\scriptstyle $}
\put(10.2,-1.15){$\scriptstyle {\rm inc}$}
\put(10.2,0.5){$\scriptstyle $}
\end{picture}
\end{center}
\vspace{0.3cm}
where $F=-\otimes_R^{\mathbf L} R_{\Ucal}$ is the derived tensor product, and
$G={\mathbf R} \Hom R{R_{\Ucal}}{-}$ is the derived Hom-functor.
}

\smallskip

\begin{proof}
By \ref{bousfield} and \ref{localiz},
$\Xcal=\Tria\Ucal$  is a smashing subcategory of $\Dcal(R)$ which gives rise
to a recollement
\vspace{0.3cm}
\begin{center}
\unitlength0.3cm
\begin{picture}(13,2)
\put(-3.7,0){{\Dcal(R)/\Xcal}}
\put(5.6,0){{\Dcal(R)}}
\put(13.3,0){{\Xcal}}
\put(3.2,-0.3){\oval(4,2.1)[b]}
\put(3.2,0.4){\oval(4,2.9)[t]}
\put(1.2,-0.3){\vector(0,1){0.3}}
\put(1.2,0.8){\vector(0,-1){0.3}}
\put(1.6,0.2){\vector(1,0){3.3}}
\put(2.7,2.3){$\scriptstyle \pi$}
\put(2.7,-1.35){$\scriptstyle $}
\put(2.7,0.5){$\scriptstyle {\rm }$}
\put(10.7,-0.3){\oval(4,2.1)[b]}
\put(10.7,0.4){\oval(4,2.9)[t]}
\put(8.7,-0.3){\vector(0,1){0.3}}
\put(8.7,0.8){\vector(0,-1){0.3}}
\put(9.1,0.2){\vector(1,0){3.3}}
\put(10.2,2.3){$\scriptstyle $}
\put(10.2,-1.15){$\scriptstyle {\rm inc}$}
\put(10.2,0.5){$\scriptstyle $}
\end{picture}
\end{center}
\vspace{0.3cm}
where  $\pi:\Dcal(R)\to\Dcal(R)/\Xcal$ is the quotient functor onto the
Verdier quotient. It is shown in \cite[5.3]{NR} that there is a (unique)
functor $T:\Dcal(R)/\Xcal\to \Dcal(R_{\Ucal})$ such that the derived tensor
product $F$ factors through $\pi$ as $F= T\circ \pi$. Moreover, combining
\cite[7.4, 6.5, 8.7 ]{NR} one obtains that
$Q^\cdot=\pi(R)$ satisfies $\Hom {\Dcal(R)/\Xcal}{Q^\cdot}{Q^\cdot}=R_{\Ucal}$
and $\Hom {\Dcal(R)/\Xcal}{Q^\cdot}{Q^\cdot[n]}=0$ for all integers $n\not=0$.
By \cite[5.6]{NR} it follows that the functor $T$ is an equivalence, so the
recollement above is equivalent to the one in the statement.\end{proof}

\medskip

\bigskip

\section{Constructing tilting objects from recollements.}\label{s1}

In this section we start with two exceptional objects coming from
the two end terms of a recollement and construct a tilting object
from them.

Recall that $\Dcal$ denotes a triangulated category with small coproducts.
Let $T_1$, $T_2$ be two exceptional objects in $\Dcal$ such that
$$(A1)\ \ \ \text{Hom}_\Dcal(T_1,T_2[k])=0\ \text{for\ all}\ k\in\mathbb{Z},$$
$$(A2)\ \ \ \text{Hom}_\Dcal(T_2,T_1[k])=0\ \ \text{for\ all}\
k\in\mathbb{Z}\setminus\{0,1\}.$$
Assumption (A2) generalizes the familiar condition on (exceptional)
modules to have projective dimension at most one.

Choose any morphism $\alpha:T_2\ra
T_1[1]$ and consider the triangle determined by $\alpha$:
$$T_1\ra T\stackrel{\gamma}{\to}T_2\stackrel{\alpha}{\to}T_1[1].$$
The next Proposition gives a necessary and sufficient condition
for when $T$ is exceptional.

\begin{prop}\label{partialtilting} With the notations above, $T$ is an
exceptional
object if and only if the homomorphism
$\text{End}_\Dcal(T_2)\oplus\text{\rm End}_\Dcal(T_1[1])\ra {\rm
Hom}_\Dcal(T_2,T_1[1])$ induced by $\alpha$, mapping $(f,g)$ to
$\alpha\circ f + g[1]\circ \alpha$, is surjective.
\end{prop}

\begin{proof} Applying $\text{Hom}_\Dcal(-,T_2[k])$ to the triangle determined
by $\alpha$ one obtains a long exact sequence
$$\ldots\ra \text{Hom}_\Dcal(T_2,T_2[k]) \ra \text{Hom}_\Dcal(T,T_2[k]) \ra
\text{Hom}_\Dcal(T_1,T_2[k]) \ra \ldots$$ By assumption
$\text{Hom}_\Dcal(T_1,T_2[k])=0$ for all integers $k$, and
$\text{Hom}_\Dcal(T_2,T_2[k])=0$ for all nonzero integers $k$.
Hence $\text{Hom}_\Dcal(T,T_2[k])=0$ for all nonzero integers $k$.
Applying $\text{Hom}_\Dcal(-,T_1[k])$ one obtains
$$\ldots \ra \text{Hom}_\Dcal(T_2,T_1[k]) \ra \text{Hom}_\Dcal(T,T_1[k]) \ra
\text{Hom}_\Dcal(T_1,T_1[k])\ra\ldots$$ By assumption
$\text{Hom}_\Dcal(T_1,T_1[k])=0$ for all $k\neq 0$, and
$\text{Hom}_\Dcal(T_2,T_1[k])=0$ for all $k\neq 0,1$. Hence
$\text{Hom}_\Dcal(T,T_1[k])=0$ for all $k\neq 0,1$.

Applying $\text{Hom}_\Dcal(T,-)$ to the triangle one obtains
$$\ldots\ra
\text{Hom}_\Dcal(T,T_1[k])\ra\text{Hom}_\Dcal(T,T[k])\ra
\text{Hom}_\Dcal(T,T_2[k])\ra\ldots$$
It follows that $\text{Hom}_\Dcal(T,T[k])=0$ for all $k\neq 0,1$,
and that $\text{Hom}_\Dcal(T,T[1])=0$ if and only if the map
$(T,\alpha):\text{Hom}_\Dcal(T,T_2)\ra \text{Hom}_\Dcal(T,T_1[1])$
induced by $\alpha$ is surjective.

Now consider the following commutative diagram
$$\begin{picture}(300,130)
\put(45,0){\footnotesize $0=\text{Hom}_\Dcal(T_1,T_2)$}
\put(191,0){\footnotesize $\text{Hom}_\Dcal(T_1,T_1[1])=0$}
\put(60,40){\footnotesize $\text{Hom}_\Dcal(T,T_2)$}
\put(191,40){\footnotesize $\text{Hom}_\Dcal(T,T_1[1])$}
\put(60,80){\footnotesize $\text{Hom}_\Dcal(T_2,T_2)$}
\put(190,80){\footnotesize $\text{Hom}_\Dcal(T_2,T_1[1])$}
\put(40,120){\footnotesize $0=\text{Hom}_\Dcal(T_1[1],T_2)$}
\put(190,120){\footnotesize $\text{Hom}_\Dcal(T_1[1],T_1[1])$}

\put(122,2){\vector(1,0){63}} \put(139,5){\footnotesize
$(T_1,\alpha)$} \put(118,42){\vector(1,0){67}}
\put(137,45){\footnotesize $(T,\alpha)$}
\put(121,82){\vector(1,0){63}} \put(137,85){\footnotesize
$(T_2,\alpha)$} \put(126,122){\vector(1,0){58}}

\put(90,115){\vector(0,-1){25}} \put(90,75){\vector(0,-1){25}}
\put(92,60){\footnotesize $\cong$} \put(90,35){\vector(0,-1){25}}
\put(220,115){\vector(0,-1){25}} \put(222,101){\footnotesize
$(\alpha,T_1[1])$} \put(220,75){\vector(0,-1){25}}
\put(222,61){\footnotesize $(\gamma,T_1[1])$}
\put(220,35){\vector(0,-1){25}}
\end{picture}$$
It is clear that $(\gamma,T_1[1])$ is surjective. Hence  $(T,\alpha)$
is surjective if and only if the morphism
$$(T_2,\alpha)\oplus(\alpha,T_1[1]):\text{\rm
End}_\Dcal(T_2)\oplus\text{End}_\Dcal(T_1[1])\ra {\rm
Hom}_\Dcal(T_2,T_1[1])$$ is surjective. \fatbox
\end{proof}

An alternative proof can be based on Lemma 2.1 in \cite{KN}.
\medskip

A morphism $\alpha:M\ra N$ in $\Dcal$ is called {\it
left-universal} if for any morphism $f:M\ra N$ there exists
$f_M:M\ra M$  such that $f=\alpha\circ f_M$, yielding the
following commutative diagram: $$\xymatrix{
  M \ar[d]_{f_M} \ar[dr]^{f}        \\
  M \ar[r]_{\alpha}  & N              }
$$
In other words, $\alpha$ is left universal if and only if the map
$\text{End}_\Dcal(M)\ra\text{Hom}_\Dcal(M,N)$ induced by $\alpha$
is surjective.

Dually one defines {\it right-universal} morphisms: $\alpha$ is
right universal if and only if the map
$\text{End}_\Dcal(N)\ra\text{Hom}_\Dcal(M,N)$ induced by $\alpha$
is surjective.

\begin{prop}\label{tiltinguniversal} Let $T_1$ and $T_2$ be two
exceptional objects in $\Dcal$ satisfying
conditions $(A1)$ and $(A2)$.
Then the following statements hold true.

$(1)$ The object $T\oplus T_2$ is exceptional if and only if the
morphism $\alpha:T_2\ra T_1[1]$ is left-universal.

$(2)$ The object $T\oplus T_1$ is exceptional if and only if the
morphism $\alpha:T_2\ra T_1[1]$ is right-universal.
\end{prop}

\begin{proof}

$(1)$ By assumption, $T_2$ has no self extensions, and as in the proof
of Proposition \ref{partialtilting} one
verifies $\text{Hom}_\Dcal(T,T_2[k])=0$ for all
$k\neq 0$. Applying $\text{Hom}_\Dcal(T_2,-)$ to the triangle
$$T_1\ra T\stackrel{\gamma}{\to}T_2\stackrel{\alpha}{\to}T_1[1]$$
one obtains a long exact sequence {\footnotesize
$$\ldots\ra\text{Hom}_\Dcal(T_2,T_2[k-1])\ra
\text{Hom}_\Dcal(T_2,T_1[k])\ra\text{Hom}_\Dcal(T_2,T[k])\ra
\text{Hom}_\Dcal(T_2,T_2[k])\ra\ldots$$}
The
assumptions $(A1)$ and $(A2)$ imply that $\text{Hom}_\Dcal(T_2,T[k])=0$
for all $k\neq 0,1$. Moreover, $\text{Hom}_\Dcal(T_2,T[1])=0$ if and
only if the map $\text{Hom}_\Dcal(T_2,T_2)\ra
\text{Hom}_\Dcal(T_2,T_1[1])$ induced by $\alpha$ is surjective,
which is equivalent to the left universality of $\alpha$. This
completes the 'only if' part. For the 'if' part notice further
that by Proposition \ref{partialtilting} the object $T$ is
exceptional if $\alpha$ is left universal.

$(2)$ follows by  similar arguments: Applying the functor Hom$(T_1,-)$ we see
that
$\text{Hom}_\Dcal(T_1,T[k])
\cong\text{Hom}_\Dcal(T_1,T_1[k])$ vanishes for
all $k\neq 0$. Next, applying $\text{Hom}_\Dcal(-,T_1[k])$
we get as in the proof of
Proposition \ref{partialtilting} that
$\text{Hom}_\Dcal(T,T_1[k])$ vanishes for all
$k\neq 0,1$. Finally, we observe that
$\text{Hom}_\Dcal(T,T_1[1])=0$ if and only if
$\alpha$ is right universal.
 \fatbox\end{proof}

\bigskip

\begin{cor}\label{cortilting} Let $T_1$, $T_2$ be
exceptional objects in $\Dcal$ satisfying $(A1)$. If
$\text{Hom}_{\Dcal}(T_2,T_1[k])=0$ for all but one integer $k=n$, then
$T_1[n]\oplus T_2$ is an exceptional object in $\Dcal$.
\end{cor}

\bigskip

Let us now assume that $\Dcal$  admits a recollement
\vspace{0.3cm}
\begin{center}
\unitlength0.3cm
\begin{picture}(13,2)
\put(-0.5,0){{\Ycal}}
\put(6,0){{\Dcal}}
\put(12.9,0){{\Xcal}}
\put(3.2,-0.3){\oval(4,2.1)[b]}
\put(3.2,0.4){\oval(4,2.9)[t]}
\put(1.2,-0.3){\vector(0,1){0.3}}
\put(1.2,0.8){\vector(0,-1){0.3}}
\put(1.6,0.2){\vector(1,0){3.3}}
\put(2.7,2.3){$\scriptstyle i^\ast$}
\put(2.7,-1.33){$\scriptstyle i^!$}
\put(2,0.5){$\scriptstyle i_\ast=i_!$}
\put(10.1,-0.3){\oval(4,2.1)[b]}
\put(10.1,0.4){\oval(4,2.9)[t]}
\put(8.1,-0.3){\vector(0,1){0.3}}
\put(8.1,0.8){\vector(0,-1){0.3}}
\put(8.4,0.2){\vector(1,0){3.3}}
\put(9.5,2.3){$\scriptstyle j_\ast$}
\put(9.5,-1.15){$\scriptstyle j_!$}
\put(9,0.5){$\scriptstyle j^!=j^\ast$}
\end{picture}
\end{center}

\bigskip

 \noindent
 Since $i_!$ and $j_!$ are full embeddings,
we identify $\mcy$ and $\mcx$ with their images under $i_!$ and
$j_!$ respectively.

\medskip

\begin{thm}\label{constructtogen}
Assume that $\Dcal$ admits a recollement as above.
Let $T_1$ be an exceptional generator of $\mcx$, and let $T_2$ be
a  tilting object in $\mcy$ such that $${\rm
Hom}_\Dcal(T_2,T_1[k]) = 0\ \text{for\ all}\ k \in\Z\setminus\{
0,1\},$$ and  $I = {\rm Hom}(T_2,T_1[1])$ is a set.
Consider the morphism $\alpha: T_2^{(I)} \ra T_1[1]$
induced by all elements of $I$, and let  $$T_1 \ra T \ra
T_2^{(I)} \stackrel{\alpha}{\ra} T_1[1]$$ be the triangle
determined by $\alpha$. Then $T \oplus T_2$ is  an exceptional
generator of $\Dcal$.
\end{thm}

\begin{proof}
First of all, note that the morphism $\alpha: T_2^{(I)} \ra
T_1[1]$ is left-universal. Indeed,  every map $f \in {\rm
Hom}_{\Dcal}(T_2,T_1[1])$ factors through $\alpha$ by
construction, and so does every map $f \in {\rm
Hom}_{\Dcal}(T_2^{(I)},T_1[1])$ by the universal property of
coproducts: $$\xymatrix{
  T_2^{(I)} \ar[d]_{} \ar[dr]^{f}        \\
  T_2^{(I)} \ar[r]_{\alpha}  & T_1              }$$

Next, we verify
that the objects $T_1$ and
$T_2^{(I)}$ in $\Dcal$ satisfy the assumptions of
Proposition \ref{tiltinguniversal} (1).
Of course, $T_1$ is an exceptional object.
Also $T_2^{(I)}$ is an exceptional object. In fact,  by the
self-compactness of $T_2$, we have ${\rm
Hom}_{\Dcal}(T_2^{(I)}, T_2^{(I)}[n]) \cong {\rm
Hom}_{\Dcal}(T_2^{(I)},T_2[n])^{(I)} \cong {{\rm
Hom}_{\Dcal}(T_2,T_2[n])^{(I)}}^{I} = 0$ for all $n \neq 0$.
Further, for all $n\in\Z$ we have $T_2^{(I)}[n]\in\Ycal$, and
we infer from the orthogonality in the recollement that
${\rm Hom}_{\Dcal}(T_1,T_2^{(I)}[n])=0$, proving  condition $(A1)$.
Condition $(A2)$ holds by  assumption, because
${\rm Hom}_{\Dcal}(T_2^{(I)},T_1[n]) \cong {\rm
Hom}_{\Dcal}(T_2,T_1[n])^I$.
Now Proposition \ref{tiltinguniversal} (1) yields that $T \oplus T_2^{(I)}$,
and thus also $T\oplus T_2$, is an exceptional object.

So, it remains to show that $T \oplus T_2$ generates $\Dcal$, or
equivalently, that $T_1\oplus
T_2$ generates $\Dcal$.
Assume that $M \in \Dcal$ satisfies ${\rm Hom}_{\Dcal}(T_1 \oplus T_2, M[n]) =
0$ for all $n$, and take the
canonical triangle defined by the  recollement of $\Dcal$
$$M_{\mcy} \ra M \ra M_{\mcx} \ra M_{\mcy}[1]$$ where $M_{\mcx}
\in \mcx$ and $M_{\mcy} \in \mcy$. Applying ${\rm Hom}_{\Dcal}(T_1,-)$
we have $${\rm Hom}_{\Dcal}(T_1, M_{\mcx}[n]) = 0\ \text{for\ all}\ n.$$
Since $\Xcal$ is generated by $T_1$, we deduce
$M_{\mcx} = 0$, whence  $M \cong M_{\mcy} \in \mcy $. Since  $\mcy$ is
generated by $T_2$, and $${\rm Hom}_{\Dcal}(T_2,M[n]) =
0\ \text{for\ all}\ n,$$  we conclude that $M = 0$.
Now the proof is complete. \fatbox
\end{proof}

\medskip

A particularly nice situation arises by adding some finiteness conditions.

\begin{thm}\label{constructto} Assume that $\Dcal$ is $K$-linear
over a field $K$. Let  $\Xcal$ be a localizing subcategory of $\Dcal$,
and  $\Ycal=\Ker\Hom{\Dcal}{\Xcal}{-}$.
 Let further $T_1, T_2\in\Dcal$ be compact objects such that $T_1$ is a tilting
object in $\mcx$, and  $T_2$ is a  tilting object in $\mcy$.
Assume that $${\rm Hom}_\Dcal(T_2,T_1[k]) = 0\ \text{for\ all}\ k
\in\Z\setminus\{ 0,1\}.$$ Furthermore, suppose that ${\rm
Hom}_{\Dcal}(T_2,T_1[1])$ is a finite dimensional $K$-vector space
 with basis $\alpha_1,\ldots,\alpha_m$. Consider
the canonical maps
$$\alpha:T_2^{\oplus m}\ra T_1[1],\ \ \ \ \beta:T_2\ra T_1[1]^{\oplus
m}$$ defined by $\alpha_1,\ldots,\alpha_m$,
and let
$$T_1\ra C_1\ra T_2^{\oplus m}\stackrel{\alpha}{\to}T_1[1]$$
$$T_1^{\oplus m}\ra C_2\ra T_2\stackrel{\beta}{\to}T_1[1]^{\oplus m}$$
be the  triangles determined by $\alpha$ and $\beta$,
respectively. Then
 $C_1\oplus T_2$ and $T_1\oplus
C_2$ are tilting objects in $\Dcal$.
\end{thm}

\begin{proof}
It is clear that $\alpha$ and $\beta$ are left and right
universal, respectively. Now the statement follows by similar
arguments as in the proof of \ref{constructtogen}. Note that here
$T_1$ and $T_2\,^{\oplus m}$ verify condition (A1) because
$T_1\in\Xcal$, $T_2\in\Ycal$, and $\Ycal$ is closed under finite
coproducts and shifts. $\fatbox$\end{proof}

This construction extends the familiar construction of a 'Bongartz
complement' \cite{Bo}.

\bigskip

\section{Some examples}

Now we apply the previous results to various situations in the literature.
In all cases, recollements come up naturally. These recollements then produce
exceptional objects or tilting objects previously constructed in different
ways. Moreover, the recollements may be used to give new proofs of some known
results; we refrain from giving details and instead just provide references.

\begin{ex}\label{S/R} Injective ring epimorphisms have been studied in
\cite{AHS} in order to construct tilting modules of projective dimension one.
We recover this construction by
 showing that the recollement induced by an injective homological  epimorphism
 produces
 the tilting object found in \cite{AHS}. \\
{\rm We have seen in \ref{ringepi} that every homological
ring epimorphism $\lambda:R\ra S$ gives rise to a recollement of
$\Dcal(R)$. Assume now that $\lambda$ is injective and that $S$ is
an $R$-module of projective dimenson at most one. Then we have a
triangle
$$S/R[-1]\to R\stackrel{\lambda}{\to} S\to$$ so the corresponding
recollement is of the form
\vspace{0.3cm}
\begin{center}
\unitlength0.3cm
\begin{picture}(14,2)
\put(-1.8,0){{$\Dcal(S)$}}
\put(5.6,0){{$\Dcal(R)$}}
\put(13.3,0){{$\Tria S/R$}}
\put(3.2,-0.3){\oval(4,2.1)[b]}
\put(3.2,0.4){\oval(4,2.9)[t]}
\put(1.2,-0.3){\vector(0,1){0.3}}
\put(1.2,0.8){\vector(0,-1){0.3}}
\put(1.6,0.2){\vector(1,0){3.3}}
\put(2.7,2.1){$\scriptstyle F$}
\put(2.7,-1.28){$\scriptstyle G$}
\put(2.7,0.5){$\scriptstyle \lambda_\ast$}
\put(10.7,-0.3){\oval(4,2.1)[b]}
\put(10.7,0.4){\oval(4,2.9)[t]}
\put(8.7,-0.3){\vector(0,1){0.3}}
\put(8.7,0.8){\vector(0,-1){0.3}}
\put(9.1,0.2){\vector(1,0){3.3}}
\put(10.2,2.3){$\scriptstyle $}
\put(10.2,-1.15){$\scriptstyle $}
\put(10.2,0.5){$\scriptstyle $}
\end{picture}
\end{center}

\bigskip



\noindent Recall from \cite{AHS} that  $S\oplus S/R$ is a tilting
$R$-module in the sense of the definition on page \pageref{tilt}.
Indeed, this is exactly the exceptional object constructed in
Corollary \ref{cortilting} from the exceptional objects  $T_1=S/R$
and $T_2=S$, since  $\text{Hom}(S,S/R[k])\neq 0$ iff $k=0$, and
$\text{Hom}(S/R,S[k])= 0$ for all $k\in\Z$; for details
cf.~\cite{AHS}. }
\end{ex}

\begin{ex} Canonical algebras are derived equivalent to categories
of coherent sheaves over weighted projective lines. In studying
these categories, homological epimorphisms play a major role, as
demonstrated by Geigle and Lenzing in \cite{GL}. We illustrate
our construction above by reviewing some results from \cite{GL}. \\
{\rm Let $A$ be a finite dimensional algebra, and $M$ a finite
dimensional right $A$-module with projective dimension $0$ or $1$.
Suppose $M$ is an exceptional module (that is, $\text{Ext}_A^1(M,M)=0$) such that
$\text{Hom}_A(M,A)=0$ and $\End(M)=K$ is a skew field. Write $m$ for
the dimension of $\text{Ext}_A^1(M,A)$ over $K$, and construct the
universal extension
$$0\ra A\ra N\ra M^{\oplus m}\ra 0.$$ Indeed, $N$ is the Bongartz
complement of $M$.

On the other hand, by assumption $M$ is a compact exceptional
object in the derived module category
$\Dcal(A)$. By \ref{bousfield} and
\ref{singleobjects}, it generates  a smashing subcategory $\Tria M$
and  a recollement of the form \vspace{0.3cm}
\begin{center}
\unitlength0.3cm
\begin{picture}(13,2)
\put(-3.2,0){{$\Tria N$}} \put(5.6,0){{$\Dcal(A)$}}
\put(13.3,0){{$\Tria M$}} \put(3.2,-0.3){\oval(4,2.1)[b]}
\put(3.2,0.4){\oval(4,2.9)[t]} \put(1.2,-0.3){\vector(0,1){0.3}}
\put(1.2,0.8){\vector(0,-1){0.3}} \put(1.6,0.2){\vector(1,0){3.3}}
\put(2.7,2.1){$\scriptstyle i^*$} 
\put(2.7,0.5){$\scriptstyle i_*$}
\put(10.7,-0.3){\oval(4,2.1)[b]} \put(10.7,0.4){\oval(4,2.9)[t]}
\put(8.7,-0.3){\vector(0,1){0.3}}
\put(8.7,0.8){\vector(0,-1){0.3}} \put(9.1,0.2){\vector(1,0){3.3}}
\put(10.2,2.3){$\scriptstyle $} \put(10.2,-1.15){$\scriptstyle $}
\put(10.2,0.5){$\scriptstyle $}
\end{picture}
\end{center}
\vspace{0.3cm} In fact, we know from  \ref{singleobjects} that
$\Ker\Hom{\Dcal}{\Tria M}{-}=\Tria i^*(A)$, and we will see in
Proposition \ref{brickreflection} in the Appendix that  $i^*(A)=N$.
In particular, $i^*(A)$ is exceptional. Hence by Proposition \ref{ringepi} the
recollement is induced by a homological ring epimorphism $\lambda:A\ra B$,
where $B=\End(N)$ is the endomorphism ring of $N$, and
$\Tria N$ is equivalent
to the derived category $\Dcal (B)$.
Since $M$ is a compact exceptional generator of $\Tria M$, we infer from
Keller's theorem in \ref{generators} that $\Tria M$ is equivalent to $\Dcal(K)$.
Thus the recollement has the form
\vspace{0.3cm}
\begin{center}
\unitlength0.3cm
\begin{picture}(13,2)
\put(-2.2,0){{$\Dcal (B)$}} \put(5.6,0){{$\Dcal(A)$}}
\put(13.3,0){{$\Dcal (K)$}} \put(3.2,-0.3){\oval(4,2.1)[b]}
\put(3.2,0.4){\oval(4,2.9)[t]} \put(1.2,-0.3){\vector(0,1){0.3}}
\put(1.2,0.8){\vector(0,-1){0.3}} \put(1.6,0.2){\vector(1,0){3.3}}
\put(2.7,2.1){$$} 
\put(2.7,0.5){$$}
\put(10.7,-0.3){\oval(4,2.1)[b]} \put(10.7,0.4){\oval(4,2.9)[t]}
\put(8.7,-0.3){\vector(0,1){0.3}}
\put(8.7,0.8){\vector(0,-1){0.3}} \put(9.1,0.2){\vector(1,0){3.3}}
\put(10.2,2.3){$\scriptstyle $} \put(10.2,-1.15){$\scriptstyle $}
\put(10.2,0.5){$\scriptstyle $}
\end{picture}
\end{center}
\vspace{0.3cm}

 We will see in Lemma \ref{reflective} and \ref{perpendicular}
that the essential image  of the restriction functor
$\lambda_\ast: \mbox{\rm Mod-$B$}\to\mbox{\rm Mod-$A$}$ coincides with the perpendicular category
$$\widehat{M}=\{ X\in\mbox{\rm Mod-$A$}\,\mid\, {\rm Hom}_A(M,X)= {\rm Ext}_A^1(M,X)=0\}$$ and that
$\lambda$  can be chosen as universal localization at $\Ucal=\{M\}$. Moreover,   $\lambda:A\ra B$, when viewed as
an $A$-module homomorphism, coincides up to isomorphism with the map $A\ra N$ in
the universal extension (this can also be deduced from the adjointness of $(i^*, i_*)$), and it is therefore injective.
By induction we recover   \cite[Theorem 4.16]{GL}.

For example, take $A$ to be a canonical algebra of weight type
$(p_1, p_2, \cdots, p_n)$, and $M$ an exceptional simple regular
module corresponding to the weight $p_i$. By
\cite[Theorem 10.3]{GL} we obtain that the algebra $B$ is Morita equivalent to the
canonical algebra of weight type $(p_1, \cdots, p_{i-1}, p_i-1,
p_{i+1}, \cdots, p_n)$.}
\end{ex}

\begin{ex} Ladkani has  constructed and
studied derived equivalences for incidence
algebras of partially ordered sets. The exceptional objects he considered in
this context \cite{L} are also produced by our construction in Section 2. \\
{\rm Let $X$ be a finite poset, $i:Y\hookrightarrow X$ a closed
subset, and $j:U\hookrightarrow X$ the open complement. Following
Ladkani's notation, we let $Sh(X)$ be the category of sheaves over
$X$ with values in the category of finite dimensional vector
spaces over a field $K$. By \cite {L} this is equivalent to the
category $\text{mod}(KX)$ of finite dimensional modules over the
incidence algebra $KX$. Let $\Dcal^b(X)=\Dcal^b(Sh(X))\cong
\Dcal^b(\text{mod}(KX))$
be the bounded derived category. 
By \cite{L}, there exists a `left' recollement of $\Dcal^b(X)$
built up by $\Dcal^b(Y)$ and $\Dcal^b(U)$ \vspace{0.3cm}
\begin{center}
\unitlength0.3cm
\begin{picture}(15,2)
\put(-3.0,0){{$\Dcal^b(Y)$}} \put(6,0){{$\Dcal^b(X)$}}
\put(14.9,0){{$\Dcal^b(U)$}} 
\put(3.2,0.4){\oval(4,2.9)[t]} 
\put(1.2,0.8){\vector(0,-1){0.3}} \put(1.6,0.2){\vector(1,0){3.3}}
\put(2.7,2.3){$\scriptstyle i^\ast$} 
\put(2,0.5){$\scriptstyle i_\ast=i_!$}
\put(12.1,-0.3){\oval(4,2.1)[b]} 
\put(10.1,-0.3){\vector(0,1){0.3}}
\put(10.4,0.2){\vector(1,0){3.3}}
\put(11.5,-1.15){$\scriptstyle j_!$} \put(11,0.5){$\scriptstyle
j^!=j^\ast$}
\end{picture}
\end{center}
\vspace{0.3cm} Take $T_2$ to be the direct sum of indecomposable
projective modules of $KY$, and $T_1$ the direct sum of
indecomposable injective modules of $KU$. One checks directly, as
in \cite[Proposition 4.5]{L}, that
$\text{Hom}(i_*(T_2),j_!(T_1)[k])\neq 0$ if and only if $k=1$.
Hence by Corollary \ref{cortilting}, $i_*(T_2)\oplus j_!(T_1)[1]$
is a tilting object in $D^b(X)$ (as  shown in \cite[Proposition 4.5]{L}). }
\end{ex}

\begin{lemma}\label{recollement}
Let $A$ be a finite dimensional algebra over
a field $K$, $e\in A$ an idempotent. Assume that the global
dimension of $eAe$ is finite, and that $Ae\lten_{eAe}eA=AeA$. Then
there exists a recollement of the form \vspace{0.3cm}
\begin{center}
\unitlength0.3cm
\begin{picture}(13,2)
\put(-5.7,0){{$\Dcal^b(A/AeA)$}} \put(6,0){{$\Dcal^b(A)$}}
\put(14.9,0){{$\Dcal^b(eAe)$}} \put(3.2,-0.3){\oval(4,2.1)[b]}
\put(3.2,0.4){\oval(4,2.9)[t]} \put(1.2,-0.3){\vector(0,1){0.3}}
\put(1.2,0.8){\vector(0,-1){0.3}} \put(1.6,0.2){\vector(1,0){3.3}}
\put(2.7,2.3){$\scriptstyle i^\ast$} \put(2.7,-1.33){$\scriptstyle
i^!$} \put(2,0.5){$\scriptstyle i_\ast=i_!$}
\put(12.1,-0.3){\oval(4,2.1)[b]} \put(12.1,0.4){\oval(4,2.9)[t]}
\put(10.1,-0.3){\vector(0,1){0.3}}
\put(10.1,0.8){\vector(0,-1){0.3}}
\put(10.4,0.2){\vector(1,0){3.3}}
\put(11.5,2.3){$\scriptstyle j_\ast$}
\put(11.5,-1.15){$\scriptstyle j_!$} \put(11,0.5){$\scriptstyle
j^!=j^\ast$}
\end{picture}
\end{center}
\vspace{0.3cm}
\end{lemma}

This follows from  \cite[Theorem 2.7 (b)]{PS}. The recollement is
the derived version of the following recollement of abelian
categories
\vspace{0.3cm}
\begin{center}
\unitlength0.3cm
\begin{picture}(13,2)
\put(-6.7,0){{$\text{mod}(A/AeA)$}} \put(6,0){{$\text{mod}(A)$}}
\put(15.9,0){{$\text{mod}(eAe)$}} \put(3.2,-0.3){\oval(4,2.1)[b]}
\put(3.2,0.4){\oval(4,2.9)[t]} \put(1.2,-0.3){\vector(0,1){0.3}}
\put(1.2,0.8){\vector(0,-1){0.3}} \put(1.6,0.2){\vector(1,0){3.3}}
\put(2.7,2.3){$\scriptstyle i^\ast$} \put(2.7,-1.33){$\scriptstyle
i^!$} \put(2,0.5){$\scriptstyle i_\ast=i_!$}
\put(13.1,-0.3){\oval(4,2.1)[b]} \put(13.1,0.4){\oval(4,2.9)[t]}
\put(11.1,-0.3){\vector(0,1){0.3}}
\put(11.1,0.8){\vector(0,-1){0.3}}
\put(11.4,0.2){\vector(1,0){3.3}}
\put(12.5,2.3){$\scriptstyle j_\ast$}
\put(12.5,-1.15){$\scriptstyle j_!$} \put(12,0.5){$\scriptstyle
j^!=j^\ast$}
\end{picture}
\end{center}
\vspace{0.3cm} where $i^*=-\otimes_A A/AeA$,
$i^!=\text{Hom}_A(A/AeA,-)$, $j_!=-\otimes_{eAe} A$,
$j^*=j^!=\text{Hom}_A(eAe,-)=-\otimes_A eAe$, and
$j_*=\text{Hom}_{eAe}(A,-)$.

\begin{ex} Let $A$ be a finite dimensional quasi-hereditary algebra and $e\in A$ a maximal
idempotent. Then the conditions in Lemma \ref{recollement} are
fullfilled, and the regular module can be
constructed from the recollement. \\
{\rm In this situation, the ideal $AeA$ generated by $e$ is a
heredity ideal. In particular it is projective as $A$-module, and
the quotient $A/AeA$ is again quasi-hereditary. Take $\tilde{T_2}$
to be the characteristic tilting module of $A/AeA$, and
$\tilde{T_1}=eAe$. Then $T_2:=i_*(\tilde{T_2})$ is the
characteristic tilting module of $A$ associated to $1-e$, and
$T_1:=j_!(\tilde{T_1})=eA$ is the projective standard module of
$A$ associated to $e$. Since $T_2$ has projective dimension at
most $1$, $\text{Hom}(T_2,T_1[k])\neq 0$ implies $k=0,1$. Consider
the right universal map $T_2\ra T_1[1]^{\oplus m}$ where
$m=\text{dimHom}(T_2,T_1[1])$ and the corresponding triangle
$$T_1^{\oplus m}\ra C_2\ra T_2\ra T_1[1]^{\oplus m}.$$
We infer from Theorem \ref{constructto} that $C_2\oplus T_1$ is a tilting
object. In fact, $C_2$ is the projective module corresponding to
$1-e$,  hence $C_2\oplus T_1$ is the regular module $A$.}
\end{ex}

\begin{ex}  Assem, Happel and Trepode \cite{AHT}
construct tilting modules for a one-point extension algebra from
tilting modules over the given algebras. We recover their construction.\\
{\rm Let $B$ be a finite dimensional algebra over an algebraically
closed field $K$, and $P_0$ a fixed projective right $B$-module
(in \cite{AHT} left modules are used). Denote by $A=B[P_0]$
the one-point extension of $B$ by $P_0$, that is, the matrix
algebra
$$A=\left[\begin{matrix} B & 0\\ P_0 & K \end{matrix}\right]$$
with ordinary matrix addition and multiplication induced from the
module structure of $P_0$. Write $e=e_B$ for the identity of $B$,
viewed as an idempotent in $A$ satisfying that $B=eAe=Ae$
and $A/AeA\cong K$. We assume the algebra $B$ has finite global
dimension. Then by Lemma \ref{recollement} there exists a
recollement of the following form \vspace{0.3cm}
\begin{center}
\unitlength0.3cm
\begin{picture}(15,2)
\put(-3.0,0){{$\Dcal^b(K)$}} \put(6,0){{$\Dcal^b(A)$}}
\put(14.9,0){{$\Dcal^b(B)$}} \put(3.2,-0.3){\oval(4,2.1)[b]}
\put(3.2,0.4){\oval(4,2.9)[t]} \put(1.2,-0.3){\vector(0,1){0.3}}
\put(1.2,0.8){\vector(0,-1){0.3}} \put(1.6,0.2){\vector(1,0){3.3}}
\put(2.7,2.3){$\scriptstyle i^\ast$} \put(2.7,-1.33){$\scriptstyle
i^!$} \put(2,0.5){$\scriptstyle i_\ast=i_!$}
\put(12.1,-0.3){\oval(4,2.1)[b]} \put(12.1,0.4){\oval(4,2.9)[t]}
\put(10.1,-0.3){\vector(0,1){0.3}}
\put(10.1,0.8){\vector(0,-1){0.3}}
\put(10.4,0.2){\vector(1,0){3.3}}
\put(11.5,2.3){$\scriptstyle j_\ast$}
\put(11.5,-1.15){$\scriptstyle j_!$} \put(11,0.5){$\scriptstyle
j^!=j^\ast$}
\end{picture}
\end{center}
\vspace{0.3cm}

Take $\tilde{T_2}=K\in\text{mod}(A/AeA)$ to be the simple module,
and $\tilde{T_1}\in \text{mod}(B)$ to be any tilting $B$-module.
Define $T_2=i_*(\tilde{T_2})$ and $T_1=j_*(\tilde{T_1})$. By the
definition of recollement $\text{Hom}_\Dcal(T_2,T_1[k])=0$ for all
$k$. Notice that $T_2$ is an injective $A$-module, hence
$\text{Hom}_\Dcal(T_1,T_2[k])=0$ for all nonzero $k$. We conclude from
 Corollary \ref{cortilting} that $T_1\oplus T_2$ is a tilting
$A$-module (as shown in \cite[ Proposition 4.1(b)]{AHT}).}
\end{ex}

\section{Constructing  recollements from tilting objects}\label{s2}

We have seen in Section \ref{s1} that recollements of the derived category can
be used to construct tilting objects or large tilting modules. We  are now
interested in the opposite direction: using tilting theory to produce
recollements. This will be achieved in the special case of tilting modules of
projective dimension one. Let us start with some preliminaries.

\subsection*{Notation.}

We fix a ring $R$ and work in the category $\M$ of all right $R$-modules.
For a class of modules $\mathcal C$ we denote
$$ {\mathcal C} ^o= \{ M \in \rmod R \mid \Hom RCM = 0 \mbox{ for all } C \in
\mathcal C  \},$$
$$ \mathcal C ^\perp= \{ M \in \rmod R \mid \Exti iRCM = 0
\mbox{ for all } C \in \mathcal C \mbox{ and all } i > 0 \}.$$ The
\emph{(right) perpendicular category} of $\mathcal C$ is denoted
by
$$\widehat{{\mathcal C}}=\mathcal C ^o\cap\mathcal C ^\perp.$$

Furthermore, we denote by
{Add}\,{$\mathcal C $}
the class consisting of all modules isomorphic to direct summands of
direct sums of modules of ${\mathcal C } $.
Finally, Gen{$\mathcal C $}
denotes the class of modules generated
 by modules of ${\mathcal C } $.

\bigskip

Recall that a subcategory $\Ycal$ of $\M$ is said to be {\em reflective} if
the inclusion ${\mathrm inc}_{\Ycal}:\Ycal\to \M$ has a left adjoint. This
means that every module $M\in\M$ admits a $\Ycal$-{\em reflection}, that is, a
morphism $\eta_M:M\to B$ such that $B\in\Ycal$ and $\Hom R{\eta_M}Y: \Hom RBY
\to \Hom RMY$
is bijective for all $Y\in\Ycal$. Of course,  $\Ycal$-reflections are
uniquely determined up to isomorphism.

\medskip

\begin{lemma}\label{reflective}
Let $\Ucal \subset \mbox{\rm Mod-$R$}$ be a  class of  modules  of
projective dimension at most one
 such that the class ${\Ucal}^\perp$ is closed under coproducts.
Then the following statements hold true.
\begin{enumerate}
\item The perpendicular category $\widehat{\Ucal}$ is closed under
products, coproducts, kernels, and cokernels. In particular,
$\widehat{\Ucal}$ is a reflective subcategory of $\M$.
\item
 There is a ring epimorphism
 $\lambda:R\to S$, which is uniquely determined up to equivalence, such that $\widehat{\Ucal}$ coincides with the essential image  of the restriction functor
$\lambda_\ast: \mbox{\rm Mod-$S$}\to\M$.
 \item  The map $\lambda:R\to S_R$, when viewed as an $R$-module homomorphism, is the  $\widehat{\Ucal}$-reflection of $R$.
 \item If $\Ucal$ consists of finitely presented modules, then $\lambda$ can be chosen as universal localization at $\Ucal$.
 \end{enumerate}
 \end{lemma}
 \begin{proof}
 (1) Clearly, $\widehat{\Ucal}$ is  closed  under direct products, and $\mathcal U^0$ is closed under  direct products and submodules, hence also under direct sums. Moreover, note that the  assumptions on $\Ucal$ imply that $\mathcal
U^{\perp}$ is a torsion class, that is, it is closed under epimorphic images and direct sums. So, we deduce  that  $\widehat{\Ucal}$ is closed under direct sums.

We now verify that  $\widehat{\Ucal}$ is closed under kernels.
Consider
$$\xymatrix{
0\ar[r]& \Ker f\ar[r]& Y\ar[rr]^{f}\ar@{->>}[dr]&& Z\\
&&& \Img f\ar@{^(->}[ur] }$$ with $Y, Z\in \widehat{\Ucal}$.
Since $\mathcal U ^{0}$ is closed under submodules and $\mathcal U
^{\perp}$ is a torsion class, we get $\Img f\in \mathcal U ^{0} \cap
\mathcal U ^{\perp}=  \widehat{\Ucal}$. Now, for $U\in \mathcal U$, applying
$\Hom RU{-}$ to the short exact
sequence $0\to{\Ker f}\to{Y}\to{\Img f}\to 0$, we get $\Ext {R}{U}{\Ker
f}=0$. This shows that  $\Ker f \in \widehat{\Ucal}$.

The closure  under cokernels is proved by similar arguments.


 Statements (2) and (3) now follow from \cite[1.2]{GP}, and statement (4) is proven in \cite[1.7]{AA}.
 \end{proof}

 \bigskip

  We now
  generalize the construction of the recollement given in Example \ref{S/R}.
Let us
  fix
 a module $M\in\M$ of projective dimension at most one
 such that $M^\perp$ is closed under coproducts.
 Set
 $$\Xcal=\Tria M$$
 and consider the orthogonal class
$$\Ycal=\Ker \Hom {\Dcal(R)}{\Xcal}{-}$$
 of all objects $Y\in\Dcal(R)$ such that
$\Hom {\Dcal(R)}{X}{Y} =0$ for all $X\in\Xcal$.

By Theorem \ref{bousfield}, the category $\Xcal$ is a localizing subcategory
of $\Dcal(R)$. Actually, it is even a smashing subcategory due to the
following observation.

 \begin{lemma}\label{perpendicular}
Let $M$ be a module of projective dimension at most one
 with corresponding stalk complex $M^\cdot$, and let ${Y^\cdot}\in\Dcal(R)$
be a complex.
The following statements are equivalent.\begin{enumerate}
\item
 $\Hom {\Dcal(R)}{M^\cdot[n]}{Y^\cdot} =0$ for all  $n\in\Z$.
\item All homologies $H^n(Y^\cdot), n\in\Z,$ belong to  $\widehat{M}$.
\end{enumerate}
\end{lemma}
\begin{proof}
Note that for a projective module $P$ and a complex
$Y^\cdot$ there is a natural isomorphism $$\Hom
{\Dcal(R)}{P}{Y^\cdot[n]} \stackrel{\sim}{\rightarrow} \Hom
{R}{P}{H^n(Y^\cdot)},\ \ \forall\ n\in\Z.$$

If the projective dimension of $M$ is zero, i.e. $M=P$ is
projective, then $\Hom {\Dcal(R)}{P}{Y^\cdot[n]}=0$ for all
$n\in\Z$ if and only if $\Hom {R}{P}{H^n(Y^\cdot)}=0$ for all
$n\in\Z$, and this is equivalent to $H^n(Y^\cdot) \in \widehat{P}$
for all $n\in\Z$.

Now suppose the projective dimension of $M$ is one. Let
$0\rightarrow P_1\stackrel{\alpha}{\rightarrow} P_0 \rightarrow
M\rightarrow 0$ be a projective resolution of $M$. Applying the
functor $\Hom {\Dcal(R)}{-}{Y^\cdot}$ to the triangle
$P_1\rightarrow P_0 \stackrel{\alpha}{\rightarrow} M \rightarrow$
 we find that $\Hom {\Dcal(R)}{M}{Y^\cdot[n]} =0$ for all
$n\in\Z$ if and only if $\Hom {\Dcal(R)}{P_0}{Y^\cdot[n]}
\stackrel{\alpha^*}{\rightarrow} \Hom {\Dcal(R)}{P_1}{Y^\cdot[n]}$ is an isomorphism
for all $n \in\Z$, if and only if $\Hom{R}{P_0}{H^n(Y^\cdot)}
\stackrel{\alpha^*}{\rightarrow} \Hom{R}{P_1}{H^n(Y^\cdot)}$ is an isomorphism for
all $n\in\Z$.  Again this is equivalent to
$H^n(Y^\cdot)\in\widehat{M}$ for all $n\in\Z$, by applying the
functor $\Hom{R}{-}{H^n(Y^\cdot)}$ to the short exact sequence
$0\rightarrow P_1\stackrel{\alpha}{\rightarrow} P_0 \rightarrow
M\rightarrow 0$.
\end{proof}

 \begin{prop}\label{Xsmashing}
Let $M\in\M$ be  a module of projective dimension at most one
 such that $M^\perp$ is closed under coproducts, and
 denote
 $\Xcal=\Tria M$. Then
 the orthogonal class
$\Ycal=\Ker \Hom {\Dcal(R)}{\Xcal}{-}$ is closed under small coproducts,
and  $\Xcal$ is   a {smashing subcategory} of $\Dcal(R)$.
\end{prop}
\begin{proof} We know from Theorem \ref{bousfield} that
$\Ycal$ is the category of all complexes $Y^\cdot$ such that \linebreak
 $\Hom {\Dcal(R)}{M^\cdot[n]}{Y^\cdot} =0$ for all  $n\in\Z$, which means by
  Lemma \ref{perpendicular}
 that all homologies $H^n(Y^\cdot), n\in\Z,$ belong to the perpendicular category $\widehat{M}$.
 Now if $(Y_i\,^\cdot)_{i\in I}$ is a family of complexes in $\Ycal$, then the $n-th$ homology of its coproduct is isomorphic to the coproduct of the $n-th$ homologies $\bigoplus_{i\in I} H^n(Y_i\,^\cdot)$ and thus belongs to $\widehat{M}$ by Lemma \ref{reflective}(3). This shows that $\Ycal $ is closed under coproducts, and thus $\Xcal$ is a smashing subcategory of $\Dcal(R)$.
 \end{proof}

\begin{cor}\label{pdimone}
Every  module $M\in\M$ of projective dimension at most one
 such that $M^\perp$ is closed under coproducts
induces
  a  {recollement}
\vspace{0.3cm}
\begin{center}
\unitlength0.3cm
\begin{picture}(20,2)
\put(-0.5,0){{\Ycal}}
\put(5.6,0){{\Dcal(R)}}
\put(13.3,0){{$\Tria M$}}
\put(3.2,-0.3){\oval(4,2.1)[b]}
\put(3.2,0.4){\oval(4,2.9)[t]}
\put(1.2,-0.3){\vector(0,1){0.3}}
\put(1.2,0.8){\vector(0,-1){0.3}}
\put(1.6,0.2){\vector(1,0){3.3}}
\put(2.7,2.3){$\scriptstyle $}
\put(2.7,-1.35){$\scriptstyle $}
\put(2.7,0.5){$\scriptstyle {\rm }$}
\put(10.7,-0.3){\oval(4,2.1)[b]}
\put(10.7,0.4){\oval(4,2.9)[t]}
\put(8.7,-0.3){\vector(0,1){0.3}}
\put(8.7,0.8){\vector(0,-1){0.3}}
\put(9.1,0.2){\vector(1,0){3.3}}
\put(10.2,2.3){$\scriptstyle $}
\put(10.2,-1.15){$\scriptstyle {\rm }$}
\put(10.2,0.5){$\scriptstyle $}
\end{picture}
\end{center}


\hfill ${\displaystyle \Box}$\end{cor}

\bigskip

\begin{ex}{\rm
Let $P$ be a finitely generated projective $R$-module. Write
$\tau_P(R)$ for the trace of $P$ in $R$ and set  $E=
\text{End}_{\Dcal(R)}\,P$. Then

(1) $P$ is a compact exceptional object, so  it induces a
recollement \vspace{0.3cm}
\begin{center}
\unitlength0.3cm
\begin{picture}(20,2)
\put(-0.5,0){{\Ycal}} \put(5.6,0){{$\Dcal(R)$}} \put(13.3,0){{$\Tria
P\sim \Dcal(E)$}} \put(3.2,-0.3){\oval(4,2.1)[b]}
\put(3.2,0.4){\oval(4,2.9)[t]} \put(1.2,-0.3){\vector(0,1){0.3}}
\put(1.2,0.8){\vector(0,-1){0.3}} \put(1.6,0.2){\vector(1,0){3.3}}
\put(2.7,2.3){$\scriptstyle $} \put(2.7,-1.35){$\scriptstyle $}
\put(2.7,0.5){$\scriptstyle {\rm }$}
\put(10.7,-0.3){\oval(4,2.1)[b]} \put(10.7,0.4){\oval(4,2.9)[t]}
\put(8.7,-0.3){\vector(0,1){0.3}}
\put(8.7,0.8){\vector(0,-1){0.3}} \put(9.1,0.2){\vector(1,0){3.3}}
\put(10.2,2.3){$\scriptstyle $} \put(10.2,-1.15){$\scriptstyle
{\rm }$} \put(10.2,0.5){$\scriptstyle $}
\end{picture}
\end{center}

\bigskip

\noindent In fact, $P$ is a tilting object in $\Tria P$. So, $\Tria P\sim \Dcal(E)$ by
Keller's  theorem in \ref{generators}.

\smallskip

(2) By Lemma \ref{reflective} the perpendicular category
$\widehat{P}$ is a reflective subcategory of $\M$. As shown in
\cite[Section 1]{CTT}, the  $\widehat{P}$ -reflection of $R$ is
$R/\tau_P(R)$, so there is a ring epimorphism $\lambda:R\ra S$ such
that $\widehat{P}$ is  the essential image of the
restriction functor $\lambda_*$, and $S_R$ as a right $R$-module
is isomorphic to $R/\tau_P(R)$.
Moreover,
$\lambda:R\to S$ can be chosen as
universal localization at $P$, or equivalently, at the zero map
$\Sigma=\{\sigma:0\to P\}$.
We can also prove this directly. Indeed, $\lambda$ is
$\Sigma$-inverting since $P\otimes_R S$ becomes zero, and it is
universal with this property, because for any $\Sigma$-inverting
ring homomorphism $\mu:R\to S'$ we have $P\otimes_R S'=0$, hence
$\tau_P(R)\otimes_R S'=0$ and therefore $\mu(\tau_P(R))=0$.

\smallskip

(3) If  $\lambda:R\ra S$ is a homological epimorphism, then by using the triangle $\tau_P(R)\to R\stackrel{\lambda}{\to} S\to$
we infer from \ref{ringepi}
  that we have a recollement
 \vspace{0.3cm}
\begin{center}
\unitlength0.3cm
\begin{picture}(17,2)
\put(-2.8,0){{$\Dcal(R_P)$}} \put(5.6,0){{$\Dcal(R)$}}
\put(13.3,0){{$\Tria \tau_P(R)$}} \put(3.2,-0.3){\oval(4,2.1)[b]}
\put(3.2,0.4){\oval(4,2.9)[t]} \put(1.2,-0.3){\vector(0,1){0.3}}
\put(1.2,0.8){\vector(0,-1){0.3}} \put(1.6,0.2){\vector(1,0){3.3}}
\put(2.7,2.3){$\scriptstyle $} \put(2.7,-1.35){$\scriptstyle $}
\put(2.7,0.5){$\scriptstyle {\rm }$}
\put(10.7,-0.3){\oval(4,2.1)[b]} \put(10.7,0.4){\oval(4,2.9)[t]}
\put(8.7,-0.3){\vector(0,1){0.3}}
\put(8.7,0.8){\vector(0,-1){0.3}} \put(9.1,0.2){\vector(1,0){3.3}}
\put(10.2,2.3){$\scriptstyle $} \put(10.2,-1.15){$\scriptstyle
{\rm }$} \put(10.2,0.5){$\scriptstyle $}
\end{picture}
\end{center}

\bigskip

This is equivalent to the recollement in (1).
Indeed, we know by Lemma \ref{perpendicular} that $\Ycal$ is the full triangulated subcategory of $\Dcal (R)$ consisting of the complexes with all homologies  in $\widehat{P}$, which is  identified with $\mbox{\rm Mod-$S$}$ by (2). The following Lemma \ref{truncation} will show that $\Ycal = \Dcal(S)$.

\bigskip

(4) If $P$ is generated by an idempotent $e\in R$, then the trace
of $P$ in $R$ is the two-sided ideal $ReR$. Hence the ring $S$ is
the quotient ring $R/ReR$, and $\lambda$ is the natural projection
$R\ra R/ReR$. Note that the latter is a homological epimorphism if
and only if $Re\lten_{eRe} eR = ReR$, and such an ideal $ReR$ is
called a {\em stratifying ideal} (see \cite[Section 2]{CPS}). In
this case we obtain a recollement \vspace{0.1cm}
\begin{center}
\unitlength0.3cm
\begin{picture}(16,2)
\put(-5.2,0){{$\Dcal(R/ReR)$}} \put(5.6,0){{$\Dcal(R)$}}
\put(13.3,0){{$\Dcal(eRe)\sim\Tria ReR$}} \put(3.2,-0.3){\oval(4,2.1)[b]}
\put(3.2,0.4){\oval(4,2.9)[t]} \put(1.2,-0.3){\vector(0,1){0.3}}
\put(1.2,0.8){\vector(0,-1){0.3}} \put(1.6,0.2){\vector(1,0){3.3}}
\put(2.7,2.3){$\scriptstyle $} \put(2.7,-1.35){$\scriptstyle $}
\put(2.7,0.5){$\scriptstyle {\rm }$}
\put(10.7,-0.3){\oval(4,2.1)[b]} \put(10.7,0.4){\oval(4,2.9)[t]}
\put(8.7,-0.3){\vector(0,1){0.3}}
\put(8.7,0.8){\vector(0,-1){0.3}} \put(9.1,0.2){\vector(1,0){3.3}}
\put(10.2,2.3){$\scriptstyle $} \put(10.2,-1.15){$\scriptstyle
{\rm }$} \put(10.2,0.5){$\scriptstyle $}
\end{picture}
\end{center}

\bigskip

This is the unbounded version of Lemma \ref{recollement}.
}
\end{ex}

\bigskip

\begin{lemma}\label{truncation}
 Let $\lambda:R\ra S$ be a homological ring
epimorphism. Then the full triangulated subcategory $\Ycal$ of
$\Dcal(R)$ consisting of those complexes whose cohomologies are
$S$-modules coincides with the essential image of the restriction
functor $\lambda_*: \Dcal (S) \ra \Dcal (R)$.
\end{lemma}

\begin{proof} We identify $\Dcal(S)$ with its image under
$\lambda_*$. It is clear that $\Dcal (S) \subset \Ycal$.
Conversely we need to show any complex in $\Ycal$ is contained in
$\Dcal(S)$. Since the restriction functor
$\lambda_*:\Dcal(S)\ra\Dcal(R)$ has both a left adjoint and a
right adjoint, the subcategory $\Dcal(S)$ of $\Dcal(R)$ is closed
under both small products and small coproducts. Therefore it is
closed under taking homotopy limits and colimits (for a definition
see the Appendix).

\medskip
By using the canonical truncation we see that any bounded complex
$M^\cdot$ is generated by its cohomology, in the sense that
$M^\cdot \in \Tria (\oplus_n H^n(M^\cdot))$. Any bounded above
complex in $\Ycal$ can be expressed as the homotopy limit of its
'quotient' complexes. These 'quotient' complexes are obtained from
the canonical truncation, and hence are bounded and generated by
their cohomologies. Since canonical truncation preserves
cohomology, the 'quotient' complexes are generated by $S_R$ in the
sense that they belong to $\Tria S_R$. Thus they belong to $\Dcal
(S)$. It follows that any bounded above complex in $\Ycal$ belongs
to $\Dcal(S)$. Dually, we express a bounded below complex in
$\Ycal$ as the homotopy colimit of its 'sub'-complexes, which are
also obtained from the canonical truncation and thus bounded and
belong to $\Tria S_R$. Since $\Tria S_R$ is closed under small
coproducts and hence closed under homotopy colimits, we see that
any bounded below complex in $\Ycal$ actually belongs to $\Tria
S_R$, which is contained in $\Dcal(S)$. Finally since any complex
is generated by a bounded above complex and a bounded below
complex by the canonical truncation, we conclude that any complex
in $\Ycal$ belongs to $\Dcal(S)$.
\end{proof}

\bigskip

Next, we consider recollements related to tilting modules.
Recall that a  module $T$ is said to be  a {\em tilting module (of projective dimension at most one)}\label{tilt}   if Gen$T=T^\perp$, or equivalently, if  the following conditions are
 satisfied:

  (T1) proj.dim$(T) \le 1$;

  (T2) $\Ext RT{T^{(I)}} = 0$ for each set $I$; and

  (T3) there is an    exact sequence
$0 \to R \to T_0 \to T_1 \to 0$ where  $T_0,T_1$ belong to Add$T$.

\smallskip

The class $T^\perp$ is then called a {\em tilting class}.
We  say that two tilting modules $T$ and $T'$ are {\em equivalent} if their tilting classes coincide.

\medskip

\begin{rem}
(1) Note that, in contrast to the definition of a tilting object,  a tilting module need not  be compact. This is the reason why one has to require the property ``exceptional'' in the stronger form of condition (T2).

\smallskip

(2) \label{ptm}
 Suppose that a module $T_1\in\M$ satisfies conditions (T1) and (T2). Then $T_1\,^\perp$ is closed under coproducts if and only if there are a set $I$ and a short  exact sequence $0\to R \to T_0\to T_1^{(I)}\to 0$ such that $T_0\oplus T_1$ is a tilting module \cite[1.8 and 1.9]{CT}.
 So, the (strongly) exceptional modules satisfying the assumptions of Corollary \ref{pdimone} are precisely the  modules $T_1$ that are  direct summands of a tilting module $T$ with $T^\perp=T_1\,^\perp$.\end{rem}

\medskip

Every tilting module is associated to a class of finitely presented modules of projective dimension  one \cite{BH} and thus to universal localization.

\begin{thm}\cite{AA} \label{aa}
For every tilting module $T$  of projective dimension one there exist
an exact sequence $$0\to R \to T_0\to T_1\to 0$$
and a set $\Ucal$ of finitely   presented modules of projective dimension  one such that
\begin{enumerate}
\item
$T_0,T_1\in\Add T$,
\item
$\Ucal^\perp=\Gen T=T_1\,^\perp$,
\item
$\widehat{\Ucal}=\widehat{T_1}$ coincides with the essential image  of the restriction functor
$ \mbox{\rm Mod-$R_{\Ucal}$}\to\M$ induced by
the universal localization \mbox{$\lambda_{\Ucal}: R\to R_{\Ucal}$.}
\end{enumerate}
\end{thm}

\medskip

We are now ready for the main result of this section. It  associates a recollement to every tilting module, and it discusses when this recollement has the properties considered in  Theorem \ref{constructtogen}.

\begin{thm}\label{recoll}
Every tilting module $T$ of projective dimension one gives rise to
 a recollement
\vspace{0.3cm}
\begin{center}
\unitlength0.3cm
\begin{picture}(13,2)
\put(0,0){{\Ycal}}
\put(5.6,0){{\Dcal(R)}}
\put(13.3,0){{\Xcal}}
\put(3.2,-0.3){\oval(4,2.1)[b]}
\put(3.2,0.4){\oval(4,2.9)[t]}
\put(1.2,-0.3){\vector(0,1){0.3}}
\put(1.2,0.8){\vector(0,-1){0.3}}
\put(1.6,0.2){\vector(1,0){3.3}}
\put(2.7,2.3){$\scriptstyle  q$}
\put(2.7,-1.35){$\scriptstyle $}
\put(2.7,0.5){$\scriptstyle {\rm inc}$}
\put(10.7,-0.3){\oval(4,2.1)[b]}
\put(10.7,0.4){\oval(4,2.9)[t]}
\put(8.7,-0.3){\vector(0,1){0.3}}
\put(8.7,0.8){\vector(0,-1){0.3}}
\put(9.1,0.2){\vector(1,0){3.3}}
\put(10.2,2.3){$\scriptstyle $}
\put(10.2,-1.15){$\scriptstyle {\rm inc}$}
\put(10.2,0.5){$\scriptstyle  $}
\end{picture}
\end{center}
\vspace{0.5cm}
with the following properties.
\begin{enumerate}

\item There is a set $\Ucal$ of finitely   presented modules of projective dimension  one
such that $\Gen T=\Ucal^\perp$, and $\Xcal=\Tria\Ucal$.

\item There is a module $T_1\in\Add T$ such that $\Gen T=T_1\,^\perp$, and  {$\Xcal=\Tria T_1$.}
In particular, $T_1$ is an exceptional generator of $\Xcal$.

\item $T_2=q(R)$ is a compact generator of $\Ycal$. Moreover, $T_2$ is a  tilting object in $\Ycal$ if and only if the universal localization $\lambda_{\Ucal}: R\to R_{\Ucal}$    of $R$ at $\Ucal$ is a homological epimorphism. In this case, there is an equivalence $ \Dcal(R_{\Ucal})\to\Ycal$, and the recollement above is equivalent to the one induced by $\lambda_{\Ucal}$.
\\
If, in addition,  the $R$-module $R_{\Ucal}$ has projective dimension at most one, then
\linebreak
$\text{\rm Hom}_{\Dcal(R)}(T_2, T_1[n])=0\ \text{\rm for\ all} \ n\not= 0,1.$

\item $T\in\m$  if and only if there are a ring $E$ and an equivalence $\mu:\Xcal\to \Dcal(E)$ such that $\mu(T_1)=E_E$. In this case, we can choose $\Ucal=\{T_1\}$.

\end{enumerate}
\end{thm}

\begin{proof}
Choose $T_1$ and $\Ucal$ as in Theorem \ref{aa}. Then $\Dcal(R)$
is a recollement  of  $\Tria T_1$ and
$\Ycal=\Ker\Hom{\Dcal(R)}{\Tria T_1}{-}$ by \ref{pdimone}, and by
\ref{bousfield} it is also a recollement $\Tria\Ucal$ and $\Ycal'
=\Ker\Hom{\Dcal(R)}{\Tria\Ucal}{-}$. Recall  that $\Ycal$ is the
category of all complexes $Y^\cdot$ such that
\linebreak
$\Hom
{\Dcal(R)}{T_1^\cdot[n]}{Y^\cdot} =0$ for all  $n\in\Z$, which
means by
  Lemma \ref{perpendicular}
 that all homologies $H^n(Y^\cdot), n\in\Z,$ belong to the perpendicular category $\widehat{T_1}$.
 Similarly, $\Ycal'$ consists of all complexes $Y^\cdot$ such that all homologies $H^n(Y^\cdot), n\in\Z,$ belong to the perpendicular category $\widehat{\Ucal}$. But $\widehat{T_1}=\widehat{\Ucal}$, thus $\Ycal=\Ycal'$, and the two recollements coincide. This proves (1) and (2).

 (3) First of all, note that the compact generator $R$ of $\Dcal(R)$ is mapped by $q$ to a compact generator $T_2^\cdot=q(R)$ of $\Ycal$, see \cite[4.3.6, 4.4.8]{N}.

 If  $T_2$ is a tilting object, then we know from \ref{ringepi} that our recollement is equivalent to the one induced by a homological ring epimorphism $\lambda:R\to S$.
 That means that $\Ycal$ coincides with the essential image of the restriction functor $\lambda_\ast:\Dcal(S)\to\Dcal(R)$. But then, using the description of $\Ycal$ given in   Lemma \ref{perpendicular},  we see that $\widehat{T_1}$ coincides with the essential image of the restriction functor $\mbox{\rm Mod-$S$}\to \M$ induced by $\lambda$. On the other hand, we know from Theorem \ref{aa} that $\widehat{T_1}=\widehat{\Ucal}$ coincides with the essential image  of the restriction functor
$ \mbox{\rm Mod-$R_{\Ucal}$}\to\M$ induced by
the universal localization $\lambda_{\Ucal}: R\to R_{\Ucal}$.
 By  the uniqueness of the ring epimorphism in Lemma \ref{reflective}(2)
 we conclude that $\lambda$ and $\lambda_\Ucal$ are in the same epiclass, and thus also
  $\lambda_\Ucal$ is a homological epimorphism.

Conversely, if $\lambda_\Ucal$ is a homological epimorphism, then we know from  \ref{univloc} that our recollement is equivalent to the one induced by $\lambda_{\Ucal}$. In particular, it follows from \ref{ringepi} that $T_2=q(R)$ is an exceptional object, hence   a compact tilting object in $\Ycal$. Moreover, $T_2$ is quasi-isomorphic to the stalk complex given by the $R$-module $R_\Ucal$.
Thus
 $\text{\rm Hom}_{\Dcal(R)}(T_2, T_1[n])\cong \Exti{n}{R}{R_\Ucal}{T_1}$ vanishes  for all $n\not= 0,1$
if pdim$R_\Ucal\le 1$.

(4) If $T\in\m$, then $T_1$ is  compact, hence a tilting object in $\Xcal$. So,
by Keller's  theorem in \ref{generators} there is a differential graded algebra $E=\mathbf{R}\Hom{}{T_1}{T_1}$ having homology concentrated in zero and $H^0(E)\cong  \text{End}_{\Dcal(R)}\,T_1$ with an equivalence $\mu:\Xcal\to \Dcal(E)$ such that $\mu(T_1)=E_E$.

Conversely, if we have an equivalence $\mu:\Xcal\to \Dcal(E)$ such that $\mu(T_1)=E_E$, then there is a fully faithful functor $\Dcal(E)\to \Dcal(R)$ mapping $E_E$ onto $T_1$. By \cite[1.7]{J} it follows that $T_1$ is compact in $\Xcal$, and we infer from \cite[4.4.8]{N} that $T_1$ is even compact in $\Dcal(R)$. But this means that $T_1$, and therefore also $T$, is in $\m$.
\end{proof}

\bigskip

\begin{rem}\label{trace} Let the assumptions and notations be as in Theorem \ref{recoll}.

(1) $R_\Ucal\cong T_0/\tau_{T_1}(T_0)$ where $\tau_{T_1}(T_0)$
denotes the trace of $T_1$ in $T_0$. This follows from
\cite[Section 1]{CTT}, since we know from Lemma
\ref{reflective}(3) that the universal localization
$\lambda_{\Ucal}: R\to R_{\Ucal}$, when viewed as an $R$-module
homomorphism, is the $\widehat{\Ucal}$-reflection of $R$.

(2) $\Ycal=\Tria T_2=\Tria R_\Ucal$. In fact, by definition the
module perpendicular category $\widehat{\Ucal}$ is a subcategory
of the triangular perpendicular category $\Ycal=\Ker \Hom
{\Dcal(R)}{\Ucal}{-}=\Tria T_2$. Hence $R_\Ucal$ belongs to $\Tria
T_2$ and $\Tria R_\Ucal\subseteq \Tria T_2$. Conversely, $\Tria
T_2$ is closed both under small coproduct by definition and under
small product since it is the right perpendicular category of
$\Tria T_1$. Using the argument of the proof of Lemma
\ref{truncation}, any complex in $\Tria T_2$ is generated by its
cohomology, and hence contained in $\Tria R_\Ucal$ by Lemma
\ref{perpendicular} and Theorem \ref{aa}.



(3) If $T_2$ is exceptional, then $R_\Ucal\cong T_2$. In fact, we see  as in the proof of the Proposition in \ref{ringepi} that  $T_2$ has homology concentrated in zero, and  yields a $\Ycal$-reflection $\eta_R:R\to T_2$.
Then $T_2\in\widehat{\Ucal}$ is also a $\widehat{\Ucal}$-reflection of $R_R$, and by the uniqueness of reflections, it must be isomorphic to $R_\Ucal$.

In \cite{AKL} we will see examples where $T_2$ and $R_\Ucal$ are not isomorphic.
 \end{rem}

\bigskip

\begin{cor}\label{zero}  With the assumptions and notations of Theorem \ref{aa}, the following statements are equivalent.
\begin{enumerate}
\item $\Hom R{T_1}{T_0}=0$.
\item There is a recollement
\begin{center}
\unitlength0.3cm
\begin{picture}(13,2)
\put(-2.8,0){{$\Tria T_0$}}
\put(5.6,0){{$\Dcal(R)$}}
\put(13.3,0){{$\Tria T_1$}}
\put(3.2,-0.3){\oval(4,2.1)[b]}
\put(3.2,0.4){\oval(4,2.9)[t]}
\put(1.2,-0.3){\vector(0,1){0.3}}
\put(1.2,0.8){\vector(0,-1){0.3}}
\put(1.6,0.2){\vector(1,0){3.3}}
\put(2.7,2.3){$\scriptstyle q$}
\put(2.7,-1.35){$\scriptstyle $}
\put(2.7,0.5){$\scriptstyle {\rm inc}$}
\put(10.7,-0.3){\oval(4,2.1)[b]}
\put(10.7,0.4){\oval(4,2.9)[t]}
\put(8.7,-0.3){\vector(0,1){0.3}}
\put(8.7,0.8){\vector(0,-1){0.3}}
\put(9.1,0.2){\vector(1,0){3.3}}
\put(10.2,2.3){$\scriptstyle j$}
\put(10.2,-1.15){$\scriptstyle {\rm inc}$}
\put(10.2,0.5){$\scriptstyle a$}

\end{picture}
\end{center}
\vspace{0.5cm}
\end{enumerate}
In this case, the recollement above is equivalent to the one induced by  $\lambda_{\Ucal}$, and $T$ is equivalent to the tilting module $R_\Ucal\oplus R_\Ucal/R$.

 \end{cor}
\begin{proof}
The implication $(2)\Rightarrow(1)$ follows from the definition of recollement.\\   $(1)\Rightarrow(2)$: We know that $\Dcal(R)$ is a recollement of $\Xcal=\Tria T_1$ and $\Ycal =\Ker\Hom{\Dcal(R)}{\Tria T_1}{-}$. Condition (1) means that $T_0\in\widehat{T_1}$, whence the stalk complex $T_0\,^\cdot$ belongs to $\Ycal$. So the exact sequence $0\to R\to T_0\to T_1\to 0$ gives rise to a triangle $$T_1\,^\cdot[-1]\to R^\cdot\to T_0\,^\cdot \to$$ where $T_1\,^\cdot[-1]\in \Xcal$ and $T_0\,^\cdot\in\Ycal$.  Now apply Theorem \ref{recoll} using that  $q(R)=T_0\,^\cdot$ is exceptional. For the last statement apply \cite[2.10]{AHS}, see also \cite[2.5]{AA}.


\end{proof}

\bigskip

\section{Examples of recollements induced by tilting modules}

We now provide some new examples of recollements, illustrating
particular features of our results and serving as counterexamples to
some questions that suggest themselves.

\begin{ex}\label{Kron} Our results go beyond
classical (finitely presented) tilting modules. Here we present two recollements
induced by   'large'
tilting modules over the Kronecker algebra.
The first one uses the Lukas tilting module; this recollement turns out to
be a hidden derived equivalence. The second example uses divisible modules;
here the resulting recollement is non-trivial, and it is not equivalent to the
 standard recollement by derived categories of vector
spaces, so it  becomes a
counterexample in the context of a Jordan-H\"older theorem for derived
categories (see \cite{AKL}).  \\
{\rm Let $R$ be  the Kronecker algebra, and
consider the preprojective component $\p$.
By the  Auslander-Reiten formula   $$ \p^\perp= {}^o\p $$ so   $\p^\perp$ is
the class of all {right} modules having no non-zero homomorphism to $\p$, or
in other words, the class of all modules that have no non-zero finitely
generated preprojective
direct summand (see \cite[Corollary 2.2]{R}).
There is an infinite dimensional tilting   module $L$ generating $\p^\perp$.
Its construction goes back to work by
Lukas, cf. \cite{L1,KT}.

The recollement  of $\Dcal(R)$ induced by $L$  is trivial. In fact, let us
take an exact sequence $$0\to R\to L_0\to L_1\to 0$$  and a set $\Ucal$ of
finitely presented indecomposable modules as in Theorem \ref{aa}, that is,
$L_0,L_1\in\Add L$,  $\Ucal^\perp=\Gen L=\p^\perp$, and $\widehat{\Ucal}=
\widehat{L_1}$. Then $\Ucal$ is contained in $ ^\perp{}(\p^\perp)$ and
therefore  in $\p$.  Observe that the indecomposable preprojective
$R$-modules, up to isomorphism, form a countable family $(P_n)_{n\in\N}$
where $P_1$ is simple projective, and each $P_n$ with $n>1$ generates all
modules having no direct summands isomorphic to one of $P_1,\ldots,P_{n-1}$,
hence in particular every module in  $\p^\perp$.  From this we deduce that
every module $X\in \p^\perp$ is generated by $\Ucal$, and thus cannot belong
to $\Ucal^o$, unless $X=0$.  Thus $\widehat{L_1}=\widehat{\Ucal}=
\Ucal^o\cap\p^\perp=0$.
But since $\Ycal$ consists of the complexes with all homologies in
$\widehat{L_1}$ by Lemma \ref{perpendicular}, this implies that $\Ycal=0$
and $\Xcal=\Dcal(R)$.

\medskip

Let us now consider the class of indecomposable regular right
$R$-modules $\tube$.
Again by the Auslander-Reiten formula the tilting class
$$\tube^\perp={}^o\tube$$
 is the torsion class of all {\em divisible   modules}, see \cite{R}.
We fix a tilting   module $W$
which generates $\tube^\perp$.
It is shown in \cite{RR} that $W$
can be chosen as the direct sum of a set of representatives of the Pr\"ufer
$R$-modules and the generic  $R$-module $G$. Moreover, there is an exact
sequence $$0\to R\to W_0\to W_1\to 0$$ where $W_0\cong G^d$, and $W_1$ is a
direct sum of  Pr\"ufer modules. Note that $\Hom R{W_1}{W_0}=0$, so we are in
the situation of Corollary \ref{zero}. Thus $W$ is equivalent to the tilting
module $R_\tube\oplus R_\tube/R$, and
there is a recollement
\vspace{0.5cm}
\begin{center}
\unitlength0.3cm
\begin{picture}(17,2)
\put(-2.5,0){{$\Dcal(R_\tube)$}}
\put(5.6,0){{\Dcal(R)}}
\put(13.3,0){{$\Tria \tube$}}
\put(3.2,-0.3){\oval(4,2.1)[b]}
\put(3.2,0.4){\oval(4,2.9)[t]}
\put(1.2,-0.3){\vector(0,1){0.3}}
\put(1.2,0.8){\vector(0,-1){0.3}}
\put(1.6,0.2){\vector(1,0){3.3}}
\put(2.7,2.3){$\scriptstyle  q$}
\put(2.7,-1.35){$\scriptstyle $}
\put(2.7,0.5){$\scriptstyle {\rm inc}$}
\put(10.7,-0.3){\oval(4,2.1)[b]}
\put(10.7,0.4){\oval(4,2.9)[t]}
\put(8.7,-0.3){\vector(0,1){0.3}}
\put(8.7,0.8){\vector(0,-1){0.3}}
\put(9.1,0.2){\vector(1,0){3.3}}
\put(10.2,2.3){$\scriptstyle $}
\put(10.2,-1.15){$\scriptstyle {\rm inc}$}
\put(10.2,0.5){$\scriptstyle  $}
\end{picture}
\end{center}
\vspace{0.5cm}
 where $R_\tube\cong \End _R W_0\cong {(\End _RG )}^{d\times d}$, see also
\cite{C1}, \cite[4.7]{AHS}.

}
\end{ex}

\bigskip

\begin{ex}\label{Matlis}
The following example shows that the recollement constructed as in Theorem \ref{recoll} from a tilting module $T$ can be induced by an injective homological epimorphism $\lambda:R\to Q$
despite the fact that $T$ is not of the form $Q\oplus Q/R$ as in \cite{AHS} or in Example \ref{S/R}.\\
{\rm Let $R$ be a commutative domain, and $Q$ its quotient field.
Denote by $\mathcal D$ the class  of all
{divisible modules}. It was shown by Facchini  \cite{facchini} that there is a
tilting
module of projective dimension one generating $\mathcal D$, namely the {Fuchs'
divisible module}
$\delta$, cf.\ \cite[\S VII.1]{FS}.
 Recall further that $\mathcal {D} = \mathcal{U}^{\perp}$ where
$${{\mathcal {U}}} = \{ R/rR \mid r \in R \}$$ denotes a set of representatives
of  all {cyclically presented modules}.
Moreover,  the exact sequence $0\to R\to \delta\to \delta/R\to 0$  has the
properties stated in Theorem \ref{aa}. In particular,  the perpendicular
category $\widehat{\delta/R}=\widehat{\Ucal}$ is the class of all divisible
torsion-free modules.

Note that the universal localization of $R$ at $\Ucal$ is given by the
injective flat epimorphism $\lambda:R\to Q$, see \cite[3.7]{AHS}. So, we
obtain a recollement of the form
\begin{center}
\unitlength0.3cm
\begin{picture}(20,2)
\put(-2,0){{\Dcal(Q)}}
\put(5.6,0){{\Dcal(R)}}
\put(13.3,0){{$\Tria\Ucal=\Tria\delta/R$}}
\put(3.2,-0.3){\oval(4,2.1)[b]}
\put(3.2,0.4){\oval(4,2.9)[t]}
\put(1.2,-0.3){\vector(0,1){0.3}}
\put(1.2,0.8){\vector(0,-1){0.3}}
\put(1.6,0.2){\vector(1,0){3.3}}
\put(2.7,2.3){$\scriptstyle  $}
\put(2.7,-1.35){$\scriptstyle $}
\put(2.7,0.5){$\scriptstyle {\rm inc}$}
\put(10.7,-0.3){\oval(4,2.1)[b]}
\put(10.7,0.4){\oval(4,2.9)[t]}
\put(8.7,-0.3){\vector(0,1){0.3}}
\put(8.7,0.8){\vector(0,-1){0.3}}
\put(9.1,0.2){\vector(1,0){3.3}}
\put(10.2,2.3){$\scriptstyle $}
\put(10.2,-1.15){$\scriptstyle {\rm inc}$}
\put(10.2,0.5){$\scriptstyle  $}
\end{picture}
\end{center}
\vspace{0.5cm}

 On the other hand, $\delta$ is not equivalent to a tilting module of the form $S\oplus S/R$ as in Example \ref{S/R},  unless $R$ is a Matlis domain, see \cite[2.11 (4)]{AHS}.
}

\end{ex}

\bigskip

\begin{ex}\label{cycle}
In the next example, we start with a tilting object, assign a recollement to it
as in Theorem \ref{recoll}, and then construct a tilting object from the
recollement as in  Theorem \ref{constructto}. The resulting tilting object is
different from the tilting object we started with. \\
{\rm Let $K$ be a field, and let $R$ be the $K$-algebra given the quiver
\begin{center}
\unitlength0.3cm
\begin{picture}(7,2)
\put(-0.4,0){{1$\bullet$}}
\put(5.6,0){{$\bullet$ 2}}
\put(3.2,-0.3){\oval(4,2.1)[b]}
\put(3.2,0.4){\oval(4,2.9)[t]}
\put(5.2,-0.3){\vector(0,1){0.3}}
\put(1.2,0.8){\vector(0,-1){0.3}}
\put(2.7,2.3){$\scriptstyle $}
\put(2.7,-1.35){$  \beta$}
\put(2.7,0.9){$  \alpha$}

\end{picture}
\end{center}
\vspace{0.5cm} with the relation $\beta\alpha=0$. Denote by $P_i,
I_i, S_i$, $i=1,2$, the indecomposable projective, injective, and
the simple right $R$-modules, and set $T=P_2\oplus S_2$. The
minimal left $\add T$-approximation of $R$ is given by the exact
sequence $$0\to R\to (P_2)^2\to S_2\to 0$$ Note that $S_2$ is the
socle of $P_2$, hence $\Hom R{S_2}{P_2} \not =0$, and $T$ is not
equivalent to a tilting module of the form $S\oplus S/R$ as in
 Example \ref{S/R}, see \cite[2.10]{AHS}.
Setting $\Ucal=\{S_2\}$, one easily verifies that $\Gen T=\Add
\{P_2, I_1, S_2\}= \Ucal^\perp$,
and
  that the perpendicular category $\widehat{\Ucal}=\Add I_1$.
  Using that the universal localization $\lambda_{\Ucal}: R\to R_{\Ucal}$, when
viewed as an $R$-module homomorphism, is the $\widehat{\Ucal}$-reflection of
$R$, we obtain
   $R_{\Ucal}\cong I_1\,^2$ as $R$-modules, and
$R_{\Ucal}\cong\End I_1\,^2\cong K^{2\times 2}$ as rings.
In particular, it follows that $\Exti{i}{R}{ R_{\Ucal}}{ R_{\Ucal}}=0$ for all
$i>0$, so $\lambda$ is a homological epimorphism by \cite[4.9]{GL}, and we
obtain a recollement of the form

\bigskip

\begin{center}
\unitlength0.3cm
\begin{picture}(13,2)
\put(-10,0){{$\Dcal(K^{2\times 2})\sim \Dcal(R_\Ucal)$}}
\put(5.6,0){{\Dcal(R)}}
\put(13.3,0){{$\Tria S_2\sim\Dcal(K)$}}
\put(3.2,-0.3){\oval(4,2.1)[b]}
\put(3.2,0.4){\oval(4,2.9)[t]}
\put(1.2,-0.3){\vector(0,1){0.3}}
\put(1.2,0.8){\vector(0,-1){0.3}}
\put(1.6,0.2){\vector(1,0){3.3}}
\put(2.7,2.3){$\scriptstyle  q$}
\put(2.7,-1.35){$\scriptstyle $}
\put(2.7,0.5){$\scriptstyle {\rm inc}$}
\put(10.7,-0.3){\oval(4,2.1)[b]}
\put(10.7,0.4){\oval(4,2.9)[t]}
\put(8.7,-0.3){\vector(0,1){0.3}}
\put(8.7,0.8){\vector(0,-1){0.3}}
\put(9.1,0.2){\vector(1,0){3.3}}
\put(10.2,2.3){$\scriptstyle $}
\put(10.2,-1.15){$\scriptstyle {\rm inc}$}
\put(10.2,0.5){$\scriptstyle  $}
\end{picture}
\end{center}
\vspace{0.5cm}
Moreover $$T_1=S_2,\ T_2=R_\Ucal\cong I_1$$
satisfy the assumptions of Theorem \ref{constructto} because $T_1$ has
injective dimension one and therefore $\Hom{\Dcal(R)}{T_2}{T_1[n]}\cong
\Exti{n}{R}{T_2}{T_1}$ vanishes for $n>1$. For $n=1$ we have
a one-dimensional space $\Hom{\Dcal(R)}{T_2}{T_1[1]}\cong
\Exti{1}{R}{I_1}{S_2}$ with basis given by the almost split sequence
$$0\to S_2\to I_2\to I_1\to 0$$ which yields the triangle $$  I_2\to I_1\to
S_2[1]\to$$
Note that applying Theorem \ref{constructto} we don't get the original
tilting module $T$, but a new tilting object, namely  $I_2\oplus I_1$.
}

\end{ex}

\bigskip

\begin{ex}
\label{nothom}
We close with an example where the universal localization $\lambda_\Ucal$ is
not a homological epimorphism.
{\rm Let $K$ be a field, and let $R$ be the  $K$-algebra given the quiver
\begin{center}
\unitlength0.3cm
\begin{picture}(12,2)
\put(-0.4,0){{1$\bullet$}}
\put(5.6,0){{$\bullet$ 2}}
\put(12,0){{$\bullet$ 3}}
\put(3.2,-0.3){\oval(4,2.1)[b]}
\put(3.2,0.4){\oval(4,2.9)[t]}
\put(5.2,-0.3){\vector(0,1){0.3}}
\put(1.2,0.8){\vector(0,-1){0.3}}
\put(9.6,-0.3){\oval(4,2.1)[b]}
\put(9.6,0.4){\oval(4,2.9)[t]}
\put(11.6,-0.3){\vector(0,1){0.3}}
\put(7.6,0.8){\vector(0,-1){0.3}}
\put(2.7,2.3){$\scriptstyle $}
\put(2.7,-1.35){$  \beta$}
\put(2.7,0.9){$  \alpha$}
\put(9.1,-1.35){$  \delta$}
\put(9.1,0.9){$  \gamma$}

\end{picture}
\end{center}
\vspace{0.5cm} with the relations
$\alpha\gamma=\delta\gamma=\delta\beta=0$ and
$\beta\alpha=\gamma\delta$. Denote by $P_i, I_i, S_i$, $i=1,2,3$,
the indecomposable projective, injective, and the simple right
$R$-modules. Indeed $P_1=\footnotesize \begin{array}{c} 1\\ 2   \\
1
\end{array}$, $P_2=\footnotesize \begin{array}{ccc}
& 2 &   \\ 1& &3 \\ &2&  \end{array}$, and
$P_3=\footnotesize\begin{array}{c}
 3   \\ 2 \end{array}$.

$R$ is quasi-hereditary with characteristic tilting module
$T'=P_1\oplus P_2\oplus S_1$.
 The minimal left $\add T'$-approximation of $R$ is given by the exact
sequence $$0\to R\to T_0\to T_1\to 0$$ where $T_0=P_1\oplus
(P_2)^2$ and $T_1=
\begin{array}{c}
 2   \\
  1
\end{array}$.
We consider the tilting module $T=T_0\oplus T_1$ and set
$\Ucal=\{T_1\}$. By Remark on page \pageref{trace}, the $R$-module
$R_\Ucal$ can be computed as $T_0/\tau_{T_1}(T_0)$ where
$\tau_{T_1}(T_0)$ denotes the trace of $T_1$ in $T_0$. It follows
that $R_\Ucal\cong S_1\oplus (P_2/S_2)^2$, which has non-trivial
self-extensions. We conclude that the universal localization  at
$\Ucal=\{T_1\}$ is not a homological epimorphism.


\smallskip

Another example for a universal localization  that  is not a homological
epimorphism  is given in \cite{NRS}, where a ring   of global dimension $\le
2$ with universal localization of global dimension $\ge 3$ is constructed.
The present example is quite different since $R_\Ucal$ is hereditary and
gldim$R=4$.

}

\end{ex}

\bigskip

\section{Appendix: Construction of the triangulated reflection}

Let $\Dcal=\Dcal(R)$ be the derived module category of a
ring $R$, and $T_1$ an exceptional object in $\Dcal$. Set
$\Ycal=\Ker \Hom {\Dcal(R)}{\Tria T_1}{-}$. We know from
\ref{localiz} and \ref{bousfield} that $\Xcal=\Tria T_1$ is a
localizing subcategory of $\Dcal$, thus the inclusion ${\rm
inc}:\Ycal\to \Dcal$ has a left adjoint $q$. We want to calculate
the $\Ycal$-reflection  of $R$, that is, a morphism $R\to q(R)$
such that  $q(R)$ lies in $\Ycal$ and the induced map
$\Hom{\Dcal}{q(R)}{Y} \rightarrow \Hom{\Dcal}{R}{Y}$ is an
isomorphism for any $Y\in \Ycal$.

\medskip
First assume the endomorphism ring of $T_1$ is a skew-field and
the extensions between $T_1$ and $R$ are free of finite rank over
the skew-field.

\medskip

\begin{prop}\label{brickreflection}
Suppose $T_1\in\Dcal$ is self-compact exceptional with
endomorphism ring $k$ being a skew-field. Suppose further that the
morphism spaces $\Hom{\Dcal}{T_1}{R[i]}$, $i\in\Z$, are finite
dimensional over $k$, and let $n_i =
\dim_k\Hom{\Dcal}{T_1}{R[i]}$. Consider the canonical map
$\alpha:S= \oplus_i T_1[-i]^{\oplus{n_i}} \rightarrow R$ given by
basis elements of these spaces. Then   the cone of $\alpha$ is a
$\Ycal$-reflection of $R$.
\end{prop}

\begin{proof}  The triangle $S
\stackrel{\alpha}{\to} R \rightarrow C \rightarrow$ gives rise to
the long exact cohomology sequence $$\ldots \rightarrow
\Hom{\Dcal}{T_1[-i]}{S} \stackrel{f}{\rightarrow}
\Hom{\Dcal}{T_1[-i]}{R} \rightarrow \Hom{\Dcal}{T_1[-i]}{C}
\rightarrow \dots$$ Here, $\Hom{\Dcal}{T_1[-i]}{S} =
\Hom{\Dcal}{T_1}{S[i]} =
\Hom{\Dcal}{T_1}{\oplus_jT_1[i-j]\,^{\oplus n_j}}$. Since $T_1$ is
self-compact and exceptional, $\Hom{\Dcal}{T_1}{\oplus_{j\neq i}
T_1[i-j]\,^{\oplus n_j}}=0$. Hence $\Hom{\Dcal}{T_1[-i]}{S}
= \Hom{\Dcal}{T_1}{T_1\,^{\oplus n_i}}$ has dimension $n_i$ over
$k$. Moreover, $\Hom{\Dcal}{T_1[-i]}{R}=\Hom{\Dcal}{T_1}{R[i]}$
also has dimension $n_i$ over $k$. By construction, $f$ is an
isomorphism. Therefore, $\Hom{\Dcal}{T_1[-i]}{C} =
\Hom{\Dcal}{T_1}{C[i]}$ vanishes for all $i$, which shows
$C\in\Ycal$.

Given $Y\in \Ycal$, apply $\Hom{\Dcal}{-}{Y}$ to the
triangle $S \rightarrow R \rightarrow C  \rightarrow$. Since
$\Hom{\Dcal}{S[i]}{Y}=0$ for all $i$, the induced map
$\Hom{\Dcal}{C}{Y} \rightarrow \Hom{\Dcal}{R}{Y}$ is an
isomorphism.
\end{proof}

In general, we have the following method.

\medskip

\begin{lem} Let $T_1\in \Dcal$ be a self-compact exceptional object and  $M\in
\Dcal$. Suppose that there is $N\in\Z$ such that
 $\Hom{\Dcal}{T_1}{M[i]}= 0$ for all $i>N$.
 Then there exists a complex $M_1\in \Dcal$ and a map $M
\ra M_1$ such that the following holds:

(i) $\Hom{\Dcal}{T_1}{M_1[i]}=0$ for all $i>N-1$;

(ii) The induced map $\Hom{\Dcal}{T_1}{M[i]}\ra
\Hom{\Dcal}{T_1}{M_1[i]}$ is an isomorphism for all $i\les N-2$,
and is injective for $i=N-1$;

(iii) The induced map $\Hom{\Dcal}{M_1}{Y} \ra \Hom{\Dcal}{M}{Y}$
is an isomorphism for all $Y\in \Ycal$.
\label{degreedescending}
\end{lem}

\begin{proof} Consider the universal triangle
$T_1[-N]^{(I)} \stackrel{\alpha}{\ra} M \rightarrow M_1
\rightarrow$ where $\alpha$ is the canonical map induced by all elements of $I=\Hom{\Dcal}{T_1}{M[N]}= 0$.
Applying $\Hom{\Dcal}{T_1}{-}$ to the
triangle, $M_1$ is seen to be as desired.
\end{proof}

\medskip

Without loss of generality we assume that $N=0$. In this way we
get a sequence of maps of complexes $M=M_0
\stackrel{\sigma_0}{\ra} M_1 \stackrel{\sigma_1}{\ra} M_2
\stackrel{\sigma_2}{\ra} \ldots \ra M_n \stackrel{\sigma_n}{\ra}
\ldots$ such that

(i) $\Hom{\Dcal}{T_1}{M_n[i]}=0$ for all $i > -n$;

(ii) The induced map $\Hom{\Dcal}{T_1[i]}{M_n}
\stackrel{(\sigma_n)_*}{\lra} \Hom{\Dcal}{T_1[i]}{M_{n+1}}$ is an
isomorphism for $i\ges n+2$, and is injective for $i=n+1$;

(iii) The induced map $\Hom{\Dcal}{M_{n+1}}{Y}
\stackrel{(\sigma_n)^*}{\lra} \Hom{\Dcal}{M_n}{Y}$ is an
isomorphism for all $Y\in \Ycal$.

\medskip

By definition, the {\em homotopy colimit} (see for example
\cite[4.4.9]{N}), here denoted by $M_{\infty}$, is given (up to
non-unique isomorphism) by the triangle
$$\oplus_{n\ges 0} M_n \stackrel{1-\sigma}{\ra} \oplus_{n\ges 0}
M_n \stackrel{\pi}{\ra} M_{\infty} \rightarrow$$ where $1-\sigma$
is defined by $(1,-\sigma_n)^{tr}$ on the $n$-th component $M_n$.
The {\em homotopy limit} is defined dually by using direct products.

\bigskip

\begin{thm} Let $T_1\in \Dcal$ be a compact exceptional object and $M\in
\Dcal$.
Suppose that there is $N\in\Z$ such that
 $\Hom{\Dcal}{T_1}{M[i]}= 0$ for all $i>N$.
  Define $M_n$ as above. Let $\iota: M_0 \ra
\oplus_{n\ges 0} M_n$ be the canonical embedding. Then $\pi \circ
\iota: M \ra M_{\infty}$ is the $\Ycal$-reflection of
$M$.
\end{thm}

\begin{proof} The colimit $M_{\infty}$ lies in $\Ycal$
iff for each integer $i$, the map
$$(1-\sigma)_*: \Hom{\Dcal}{T_1[i]}{\oplus_{n\ges 0} M_n} \ra
\Hom{\Dcal}{T_1[i]}{\oplus_{n\ges 0} M_n}$$ is bijective. If $i<0$
then by construction $\Hom{\Dcal}{T_1[i]}{M_n}=0$ for all $n\ges
0$, and hence $\Hom{\Dcal}{T_1[i]}{\oplus_{n\ges 0} M_n}$ is zero
as $T$ is compact. Now assume $i \ges 0$. It follows from the
construction that $\Hom{\Dcal}{T_1[i]}{M_n}=0$ for all $n>i$.
Hence $\Hom{\Dcal}{T_1[i]}{\oplus_{n> i} M_n}=0$ and
$$\Hom{\Dcal}{T_1[i]}{\oplus_{n\ges 0} M_n} =
\Hom{\Dcal}{T_1[i]}{\oplus_{n=0}^{i} M_n} = \oplus_{n=0}^{i}
\Hom{\Dcal}{T_1[i]}{M_n}.$$ The map $(1-\sigma)_*$ is given by
$$(1-\sigma)_*(f_0, f_1, \ldots, f_i) = (f_0, f_1-\sigma_0\circ f_0, \ldots,
f_i-\sigma_{i-1}\circ f_{i-1}).$$ It is straightforward now to see
the bijectivity .

\medskip

It remains to prove that $\pi \circ \iota: M \ra M_{\infty}$ is
the reflection of $M$, that is, for any $Y \in \Ycal$, the
induced map $(\pi\circ\iota)^*: \Hom{\Dcal}{M_{\infty}}{Y} \ra
\Hom{\Dcal}{M}{Y}$, sending $f$ to $f\circ\pi\circ\iota$, is
bijective.

\smallskip
Take any map $g_0: M=M_0\ra Y$. By the construction of $M_n$,
there exists uniquely for each $n\ges 0$ a map $g_n: M_n \ra Y$
such that $g_{n-1} = g_n \circ \sigma_{n-1}$. Write $g$ for the
map $(g_n)_n: \oplus_{n\ges 0} M_n \ra Y$. Then $g_0 = \iota^*(g)
= g \circ \iota$.
Apply $\Hom {\Dcal}{-}{Y}$ to the triangle $$\oplus_{n\ges 0} M_n
\stackrel{1-\sigma}{\ra} \oplus_{n\ges 0} M_n \stackrel{\pi}{\ra}
M_{\infty} \rightarrow$$ to obtain a long exact sequence
$$\ldots\ra \Hom {\Dcal}{\oplus_{n\ges 0} M_n}{Y[-1]}
\stackrel{(1-\sigma)^*}{\lra} \Hom {\Dcal}{\oplus_{n\ges
0}M_n}{Y[-1]} \ra \Hom{\Dcal}{M_{\infty}}{Y}
\stackrel{\pi^*}{\ra}$$ $$\Hom {\Dcal}{\oplus_{n\ges 0} M_n}{Y}
\stackrel{(1-\sigma)^*}{\lra} \Hom {\Dcal}{\oplus_{n\ges 0}M_n}{Y}
\ra \ldots$$ where $(1-\sigma)^* (h_n)_n =
(h_n-h_{n+1}\circ\sigma_n)_n$. It is clear that the map $g$
constructed above lies in the kernel of $(1-\sigma)^*$, and
$\iota^*: \Hom {\Dcal}{\oplus_{n\ges 0}M_n}{Y} \ra
\Hom{\Dcal}{M}{Y}$, when restricted on $\Ker (1-\sigma)^*={\rm Im}\, \pi^*$,
becomes a bijection onto $\Hom{\Dcal}{M}{Y}$.

\smallskip

Note that the map $(1-\sigma)^*: \Hom {\Dcal}{\oplus_{n\ges 0}
M_n}{Y[-1]} \ra \Hom {\Dcal}{\oplus_{n\ges 0}M_n}{Y[-1]}$ is
surjective, because all
$(\sigma_n)^*:\Hom{\Dcal}{M_{n+1}}{Y[-1]} \ra
\Hom{\Dcal}{M_n}{Y[-1]}$,   $n\ges 0$, are isomorphisms. Hence $\pi^*:
\Hom{\Dcal}{M_{\infty}}{Y} \ra \Hom{\Dcal}{\oplus_{n\ges
0}M_n}{Y}$ is an injection. Combining this with  the arguments above, we
obtain the bijectivity of $(\pi\circ\iota)^*:
\Hom{\Dcal}{M_{\infty}}{Y} \ra \Hom{\Dcal}{M}{Y}$.
\end{proof}


\medskip

In the situation of Theorem \ref{recoll}, this method can be used
for computing $q(R)$, the $\Ycal$-reflection of $R$. In
particular, it follows immediately that $q(R)$ is right bounded,
namely it belongs to $\Dcal^-(R)$.


\end{document}